\documentclass[final,leqno,showlabe]{siamltex}

\input psfig.sty

\usepackage{cite}
\usepackage{amsmath}
\usepackage{graphicx,wrapfig}
\usepackage{amssymb,amsgen,color}
\usepackage{amsfonts}
\usepackage{hyperref}
\usepackage{multirow}
\usepackage{epsfig,epstopdf,color,bm}
\usepackage{booktabs}
\usepackage{verbatim}
\usepackage{multicol}
\usepackage{makecell}

\renewcommand{\theequation}{\arabic{section}.\arabic{equation}}
\renewcommand{\theequation}{\thesection.\arabic{equation}}
\renewcommand {\thefigure} {\thesection.\arabic{figure}}

\newcommand{\bx}{\textbf{x}}
\newcommand{\nn}{\nonumber}
\newcommand{\eps}{\varepsilon}
\newtheorem{remark}{Remark}[section]
\graphicspath{{figures/}}

\title{Improved uniform error bounds on time-splitting methods for long-time dynamics of the nonlinear Klein--Gordon equation with weak nonlinearity\thanks{
{The work of the first and third authors was partially supported by Ministry of Education of Singapore grant R-146-000-290-114. The work of the second author was partially supported by NSFC grant 11771036.}}}

\author{Weizhu Bao\thanks{Department of Mathematics,
National University of Singapore, Singapore 119076 
  (\tt{matbaowz@nus.edu.sg}, \it{http://blog.nus.edu.sg/matbwz/}).}
\and Yongyong Cai\thanks{Laboratory of Mathematics and Complex Systems (Ministry of Education), School of Mathematical Sciences, Beijing Normal University, Beijing 100875, People's Republic of China.  (\tt{yongyong.cai@bnu.edu.cn}).}
\and Yue Feng\thanks{Department of Mathematics,
National University of Singapore, Singapore 119076 
  (\tt{fengyue@u.nus.edu}).}}

\begin{document}

\maketitle
\begin{abstract}
We establish improved uniform error bounds on time-splitting methods for the long-time dynamics of the nonlinear Klein--Gordon equation (NKGE) with weak cubic nonlinearity, whose strength is characterized by $\varepsilon^2$ with  $0 < \varepsilon \leq 1$ a dimensionless parameter. Actually, when $0 < \varepsilon \ll 1$, the NKGE with $O(\eps^2)$ nonlinearity and $O(1)$ initial data is equivalent to that with $O(1)$ nonlinearity and small initial data of which the amplitude is at $O(\varepsilon)$. We begin with a semi-discretization of the NKGE by the second-order time-splitting method, and followed by a full-discretization via the Fourier spectral method in space. Employing the regularity compensation oscillation (RCO) technique which controls the high frequency modes by the regularity of the exact solution and analyzes the low frequency modes by phase cancellation and energy method, we carry out the improved uniform error bounds at $O(\eps^2\tau^2)$ and $O(h^m+\eps^2\tau^2)$ for the second-order semi-discretization and full-discretization up to the long time $T_\eps = T/\varepsilon^2$ with $T$ fixed, respectively. Extensions to higher order time-splitting methods and the case of an oscillatory complex NKGE are also discussed. Finally, numerical results are provided to confirm the improved error bounds and to demonstrate that they are sharp.
\end{abstract}

\begin{keywords}
nonlinear Klein--Gordon equation, long-time dynamics, time-splitting methods, improved uniform error bounds, regularity compensation oscillation (RCO)
\end{keywords}

\begin{AMS}
35L70, 65M12, 65M15, 65M70, 81-08
\end{AMS}

\section{Introduction}
In this paper, we consider the following nonlinear Klein--Gordon equation (NKGE) \cite{CE,HN,KS,LH,Sun}
\begin{equation}
\label{eq:WNE}
\begin{cases}
\partial_{tt}u(\bx, t)-\Delta u({\bx}, t)+ u(\bx, t)+\varepsilon^2 u^3({\bx}, t)=0,&\bx \in \Omega,\quad t > 0,\\
u(\bx, 0) = u_0(\bx), \quad \partial_t u(\bx, 0) = u_1(\bx),&{\bx} \in \Omega.
\end{cases}
\end{equation}
Here, $t$ is time, $\bx$ is the spatial coordinate, $\Delta$ is the Laplace operator, $u:= u(\bx, t)$ is a real-valued scalar field, $\varepsilon\in (0, 1]$ is a dimensionless parameter used to characterize the nonlinearity strength and $\Omega =  \prod_{i = 1}^{d} (a_i, b_i) \subset \mathbb{R}^d$ $(d = 1, 2, 3)$ is a bounded domain equipped with periodic boundary conditions. The initial data $u_0(\bx)$ and $u_1(\bx)$ are two given real-valued functions independent of $\varepsilon$.

When $0 < \eps \ll 1$, by introducing $w(\bx, t)=\varepsilon u(\bx, t)$, the NKGE \eqref{eq:WNE} with weak nonlinearity and $O(1)$ initial data could be reformulated into the following NKGE with small initial data and $O(1)$ nonlinearity as
\begin{equation}\label{eq:SIE}
\begin{cases}
\partial_{tt} w({\bx}, t)-\Delta w({\bx}, t)+ w({\bx}, t)+w^3({\bx}, t) = 0, &\bx \in \Omega,\quad t > 0, \\
w({\bx}, 0) = \varepsilon u_0({\bx}),\quad \partial_t w({\bx}, 0) = \varepsilon u_1({\bx}),& {\bx} \in \Omega.
\end{cases}
\end{equation}
In fact, the long-time dynamics of the NKGE \eqref {eq:SIE} with small initial data and $O(1)$ nonlinearity is equivalent to that of the NKGE \eqref{eq:WNE} with weak nonlinearity and $O(1)$ initial data.

The nonlinear Klein--Gordon equation as a fundamental physical equation describing the motion of the spinless particle has been extensively investigated from both analytical and numerical perspectives \cite{BCZ,BD,BZ,CCLM,DSA,DB,LV,SJJ,SV,Tao}. Recently, the long-time dynamics of the NKGE \eqref{eq:WNE} in the weak nonlinearity strength regime (or \eqref{eq:SIE} with small initial data) have attracted much attention. According to the analytical results, the life-span of a smooth solution to the NKGE \eqref{eq:WNE} (or \eqref{eq:SIE}) is at least up to the time at $O(\varepsilon^{-2})$ \cite{BFG,D,D2,DS,FZ,LH}. For the long-time dynamics, near-conservation (or approximate preservation) of energy, momentum and harmonic actions has been established for the semi-discretization and  full-discretization of the NKGE \eqref{eq:SIE} with small initial data via the technique of modulated Fourier expansions \cite{CHL,CHL2,HL}. In our recent work, long-time error bounds have been rigorously established for the finite difference time domain (FDTD) methods \cite{BFY,Feng}, the exponential wave integrator Fourier pseudospectral (EWI-FP) method \cite{FY} and the time-splitting Fourier pseudospectral (TSFP) method \cite{BaoFS}. In the numerical simulations, we surprisingly found the improved uniform error bounds for the TSFP method which are better than the analytical results \cite{BaoFS}. For the long-time dynamics of the Schr\"odinger/nonlinear Schr\"odinger equation, a new technique of the {\bf regularity compensation oscillation} (RCO) has been introduced to establish the improved uniform error bounds for the TSFP method in the long-time regime \cite{BCF}. The aim of this paper is to analyze the errors of the time-splitting methods carefully and carry out the improved uniform error bounds on the semi-discretization and full-discretization for the long-time dynamics of the NKGE with the help of the RCO technique. For the refined analysis, we first reformulate the NKGE into a relativistic nonlinear Schr\"odinger equation (NLSE). According to the RCO approach, we choose a frequency cut-off parameter $\tau_0$ and control the high frequency modes ($> 1/\tau_0$) by the smoothness of the exact solution and analyze the low frequency modes ($\leq 1/\tau_0$) by phase cancellation and  energy method.

The rest of the paper is organized as follows. In section 2, we adopt the time-splitting method to discretize the NKGE in time and establish the improved uniform error bounds for the semi-discretization up to the time at $O(1/\eps^2)$. In section 3, the full-discretization by the Fourier spectral method in space is shown with the proof of improved uniform error bounds. Extensions to the complex NKGE with a general power nonlinearity and an oscillatory complex NKGE are presented in section 4. Numerical results for the long-time dynamics and the oscillatory complex NKGE are shown in section 5. Finally, some conclusions are drawn in section 6. Throughout this paper, the notation $A \lesssim B$ is used to represent that there exists a generic constant $C>0$ independent of the mesh size $h$, time step $\tau$, $\varepsilon$ and $\tau_0$ such that $|A| \leq C B$.
\section{Semi-discretization and improved uniform error bounds}
In this section, we utilize the time-splitting method to discretize the NKGE \eqref{eq:WNE} in time and establish the improved uniform error bounds up to the time at $O(1/\eps^2)$. For the simplicity of presentation, we only present the numerical schemes and corresponding results in one dimension (1D). Generalization to higher dimensions is straightforward and results remain valid without modifications. In 1D, the NKGE \eqref{eq:WNE} with periodic boundary conditions on the domain $\Omega = (a, b)$ collapses to
\begin{equation}
\label{eq:21}
\begin{cases}
\partial_{tt} u(x, t) - \partial_{xx} u(x, t)+  u(x, t) + \varepsilon^2 u^3 (x, t)= 0,& a < x < b,\ t > 0, \\
u(x, 0) = u_0(x), \ \partial_t u(x, 0) =u_1(x) , &x \in \overline{\Omega} = [a, b],
\end{cases}
\end{equation}
with boundary conditions as $u(a, t) = u(b, t), \ \partial_x u(a, t)=\partial_x u(b, t)$ for $t>0$.

For an integer $m\ge 0$, we denote  $H^m(\Omega)$ as  the set of functions $u(x)\in L^2(\Omega)$ with finite $H^m$-norm $\|\cdot\|_m$ given by
\begin{equation}
\label{sn}
\|u\|_m^2=\sum\limits_{l \in \mathbb{Z}} (1+\mu_l^2)^m|\widehat{u}_l|^2,\quad \mathrm{for}\quad u(x)=\sum\limits_{l\in \mathbb{Z}} \widehat{u}_l e^{i\mu_l(x-a)},\quad \mu_l=\frac{2\pi l}{b - a},
\end{equation}
where $\widehat{u}_l (l\in \mathbb{Z})$ are the Fourier  coefficients  of the function $u(x)$ \cite{BCZ,BaoFS}. In fact, the $H^m(\Omega)$  space is the subspace of  classical Sobolev space $W^{m,2}(\Omega)$, which consists of functions with derivatives of order up to $m-1$ being $(b - a)$-periodic. Since we consider the periodic boundary conditions, the above space $H^m(\Omega)$ is suitable.
In  addition, the space is  $L^2(\Omega)$  for $m=0$  and the corresponding norm is denoted as $\|\cdot\|$. Here, the space $H^s(\Omega)$ with $s\in\mathbb{R}$ is also well-defined consisting of functions with finite norm $\|\cdot\|_s$ \cite{STL}.


Denote $X_N := \{u = (u_0, u_1, \ldots, u_N)^T \in \mathbb{C}^{N+1} \ | \ u_0 = u_N\}$, $C_{\rm per}(\Omega)=\{u \in C(\overline \Omega) \ |\ u(a) = u(b)\}$ and
\[
Y_N := \text{span}\left\{e^{i\mu_l(x-a)},\ x \in \overline{\Omega}, \ l \in \mathcal{T}_N\right\},\ \mathcal{T}_N = \left\{l ~|~ l = -\frac{N}{2}, -\frac{N}{2}+1, \ldots, \frac{N}{2}-1\right\}.
\]
For any $u(x) \in C_{\rm per}(\Omega)$ and a vector $u \in X_N$, let $P_N: L^2(\Omega) \to Y_N$ be the standard $L^2$-projection operator onto $Y_N$, $I_N : C_{\rm per}(\Omega) \to Y_N$ or $I_N : X_N \to Y_N$ be the trigonometric interpolation operator \cite{STL}, i.e.,

\begin{equation}
P_N u(x) = \sum_{l \in \mathcal{T}_N} \widehat{u}_l e^{i\mu_l(x-a)},\qquad I_N u(x) = \sum_{l \in \mathcal{T}_N} \widetilde{u}_l e^{i\mu_l(x-a)},\qquad x \in \overline{\Omega},
\end{equation}
where
\begin{equation}
\widehat{u}_l = \frac{1}{b - a}\int^{b}_{a} u(x) e^{-i\mu_l (x-a)} dx, \quad \widetilde{u}_l = \frac{1}{N}\sum_{j=0}^{N-1} u_j e^{-i\mu_l (x_j-a)}, \quad l \in \mathcal{T}_N,
\end{equation}
with $u_j$ interpreted as $u(x_j)$ when involved.

Define the operator $\langle \nabla \rangle=\sqrt{1-\Delta}$ through its action in the Fourier space by \cite{BaoFS,BFS,FS}:
\begin{equation*}
\langle \nabla \rangle u(x)=\sum\limits_{l\in\mathbb{Z}}\sqrt{1+\mu_l^2}\widehat{u}_l e^{i\mu_l(x-a)}, \quad \mathrm{for}\quad u(x)=\sum\limits_{l\in\mathbb{Z}} \widehat{u}_l e^{i\mu_l(x-a)},\quad x\in \overline{\Omega},
\end{equation*}
and the inverse operator $ \langle \nabla \rangle^{-1}$  as
$\langle \nabla \rangle^{-1} u(x)=\sum\limits_{l\in\mathbb{Z}}\frac{\widehat{u}_l}{\sqrt{1+\mu_l^2}} e^{i\mu_l(x-a)}$,
which leads to $\| \langle \nabla \rangle^{-1} u\|_{s} =\|u\|_{s-1}\le \|u\|_{s}$.

Introduce $v(x, t) = \partial_t u(x, t)$ and
\begin{equation}
\psi(x, t) = u(x, t) - i\langle\nabla\rangle^{-1}v(x, t),\quad x\in \overline{\Omega}, \quad t\ge 0,
\label{eq:psi}
\end{equation}
then the NKGE \eqref{eq:21} could be reformulated into the following relativistic NLSE
for  $\psi:=\psi(t)=\psi(x,t)$ (spatial variable $x$ may be omitted for brevity) as
\begin{equation}
\label{eq:NLSE}
\left\{
\begin{aligned}
&i\partial_t \psi(x, t) + \langle\nabla \rangle \psi(x, t) + \frac{\eps^2}{8}\langle\nabla \rangle^{-1}  \left(\psi + \overline{\psi}\right)^3(x, t) = 0,\quad x\in\Omega, \, t > 0,\\
&\psi(a, t)=\psi(b, t), \quad \partial_x \psi(a, t)=\partial_x \psi(b, t), \quad
t \ge 0,\\
&\psi(x,0)=\psi_0(x):= u_0(x) - i\langle\nabla \rangle^{-1}u_1(x), \quad x\in \overline{\Omega},
\end{aligned}\right.
\end{equation}
where $\overline{\psi}$ denotes the complex conjugate of $\psi$. According to \eqref{eq:psi},  the solution of the NKGE \eqref{eq:21} could be recovered by
\begin{equation}
\label{uxtpsit}
u(x,t)=\frac{1}{2}\left(\psi(x,t)+\overline{\psi}(x,t)\right), \qquad
v(x, t)=\frac{i}{2}\langle\nabla \rangle\left(\psi(x,t)-\overline{\psi}(x,t)\right).
\end{equation}

\subsection{The time-splitting method}
By the splitting technique \cite{LU,MQ}, the relativistic NLSE \eqref{eq:NLSE} is split to the linear part and nonlinear part. The evolution operator for the linear part $\partial_t\psi(x,t)=i\langle\nabla\rangle\psi(x,t)$ with initial data $\psi(x,0)=\psi_0(x)$ is given by
\begin{equation}\label{eq:step1}
\psi(\cdot, t) = \varphi^{t}_T(\psi_0) := e^{it\langle\nabla \rangle}	\psi_0, \quad t \geq 0,
\end{equation}
and the nonlinear part $\partial_t\psi(x,t)=F(\psi(x,t))$ with initial data $\psi(x,0)=\psi_0(x)$ can be integrated exactly in time as
\begin{equation}\label{eq:step2}
\psi(x, t) = \varphi_V^t(\psi_0) := \psi_0 + \eps^2tF(\psi_0), \quad t \geq 0	,
\end{equation}
where the nonlinear operator $F$ is given by
\begin{equation}\label{eq:F}
F(\phi)= i\langle\nabla\rangle^{-1}G(\phi),\quad G(\phi)=\frac{1}{8}\left(\phi+\overline{\phi}\right)^3.
\end{equation}

Let $\tau > 0$  be the time step size and $t_n=n\tau$ ($n=0, 1,\ldots$) as the time steps. Denote $\psi^{[n]}:= \psi^{[n]}(x)$ as the  approximation of $\psi(x, t_n)$, then the second-order discrete-in-time splitting method via the Strang splitting for the relativistic NLSE \eqref{eq:NLSE} could be written as \cite{Strang}
\begin{equation}
\label{S2}
\psi^{[n+1]}=\mathcal{S}_{\tau}(\psi^{[n]})= \varphi^{\frac{\tau}{2}}_T \circ \varphi^{\tau}_V \circ \varphi^{\frac{\tau}{2}}_T (\psi^{[n]})=e^{i\tau\langle\nabla \rangle}\psi^{[n]}+\eps^2\tau e^{i\frac{\tau\langle\nabla \rangle}{2}}	F(e^{i\frac{\tau\langle\nabla \rangle}{2}}\psi^{[n]}),
\end{equation}
with $\psi^{[0]}=\psi_0=u_0-i\langle\nabla\rangle^{-1}u_1$. Noticing \eqref{uxtpsit}, the semi-discretization of the NKGE \eqref{eq:21} is given by
\begin{equation}
\label{uxtpsita}
u^{[n]}=\frac{1}{2}\left(\psi^{[n]}+\overline{\psi^{[n]}}\right), \quad
v^{[n]} = \frac{i}{2}\langle\nabla \rangle\left(\psi^{[n]}-\overline{\psi^{[n]}}\right),\quad n=0, 1,\ldots,
\end{equation}
where $u^{[n]}:= u^{[n]}(x)$ and $v^{[n]} := v^{[n]}$ are the approximations of $u(x, t_n)$ and $\partial_t u(x, t_n)$, respectively.
\begin{remark}
\label{remark:split}
The split-steps \eqref{eq:step1} and \eqref{eq:step2} are equivalent to the splitting of NKGE \eqref{eq:21} (in terms of $u$ and $v=\partial_tu$), respectively as
\begin{equation}
\partial_t\begin{bmatrix}u\\ v\end{bmatrix}
=\begin{bmatrix}0&1\\ \partial_{xx}-1&0\end{bmatrix}\begin{bmatrix}u\\ v\end{bmatrix},\quad
\partial_t\begin{bmatrix}u\\ v\end{bmatrix}=\begin{bmatrix}0&0\\ -\eps^2 u^2&0\end{bmatrix}\begin{bmatrix}u\\ v\end{bmatrix}.
\end{equation}
\end{remark}
\subsection{Improved uniform error bounds}
By discussions in \cite{BaoFS,D,FZ} and references therein, we make the following assumptions on the exact solution $u:=u(x, t)$ of the NKGE \eqref{eq:21} up to the time at $T_{\varepsilon} = T/\varepsilon^2$ with $T>0$ fixed:
\begin{equation*}
{\rm(A)}\qquad
\|u\|_{L^{\infty}\left([0, T_{\varepsilon}]; H^{m+1}\right)} \lesssim 1,\quad \;\; \|\partial_t u\|_{L^{\infty}\left([0, T_{\varepsilon}]; H^{m}\right)} \lesssim 1,\quad m\ge1.
\end{equation*}
 Let $u^{[n]}$ and $v^{[n]}$ be the numerical approximations obtained from the Strang splitting method \eqref{S2} with \eqref{uxtpsita}. According to the analysis in \cite{BaoFS}, under the assumption (A),  for sufficiently small $0<\tau\leq\tau_c$ ($\tau_c$ is a constant), there exists a constant $M>0$ depending on $T$, $\|u_0\|_{m+1}$, $\|u_1\|_m$, $\|u\|_{L^{\infty}([0, T_\eps]; H^m)}$ and $\|\partial_t u\|_{L^{\infty}([0, T_\eps]; H^m)}$ such that
\begin{equation}
\label{reg_semi}
\|u^{[n]}\|_{m + 1}^2 + \|v^{[n]}\|_{m}^2 \le M, \text{or equivalently}\,\|\psi^{[n]}\|_{m + 1}^2  \le M, \quad 0 \leq n \leq \frac{T/\varepsilon^2}{\tau}.
\end{equation}
The main result of this work is to establish the following improved uniform error bounds for the Strang splitting method up to the long time $T_\eps$.
\begin{theorem}
\label{thm:semi}Under the  assumption (A),  for $0 < \tau_0 \leq 1$ sufficiently small and independent of $\varepsilon$ such that,
 when $0< \tau < \alpha \frac{\pi (b-a)\tau_0}{2\sqrt{\tau^2_0(b-a)^2 + 4\pi^2(1+\tau_0^2)}}$ for a fixed constant $\alpha \in (0, 1)$, we have the following improved uniform error bounds
\begin{equation}
\|u(\cdot, t_n) - u^{[n]}\|_1 + \|\partial_t u(\cdot, t_n) - v^{[n]}\| \lesssim  \varepsilon^2 \tau^2 + \tau_0^{m+1},\quad 0 \leq n \leq \frac{T/\varepsilon^2}{\tau}.
\label{eq:error_semi}
\end{equation}
In particular, if the exact solution is sufficiently smooth, e.g. $u, \partial_t u \in H^{\infty}$, the last term $\tau_0^{m+1}$ decays exponentially fast  ($\sim e^{-c/\tau_0}$) and could be ignored practically for small enough $\tau_0$ , where the improved uniform error bounds for sufficiently small $\tau$ could be
\begin{equation}
\|u(\cdot, t_n) - u^{[n]}\|_1 + \|\partial_t u(\cdot, t_n) - v^{[n]}\| \lesssim  \varepsilon^2 \tau^{2},\quad 0 \leq n \leq \frac{T/\varepsilon^2}{\tau}.
\label{eq:inf}
\end{equation}
\end{theorem}
\begin{remark} $\tau_0\in(0,1)$ is a parameter introduced in analysis and the requirement on $\tau$ (essentially $\tau\lesssim\tau_0$) enables the improved estimates on the low Fourier modes $|l|\leq1/\tau_0$, where the constant in front of $\varepsilon^2\tau^2$ depend on $\alpha$. $\tau_0$ can be arbitrary as long as the assumed relation between $\tau$ and $\tau_0$ holds, i.e. $\tau_0$ could be fixed, or depending on $\tau$,  e.g.  $\tau_0=\frac{2\sqrt{8\pi^2+(b-a)^2}}{\alpha(b-a)\pi}\tau$.
\end{remark}
\begin{remark} Compared to the previous uniform estimates  $\|u(\cdot, t_n) - u^{[n]}\|_1 + \|\partial_t u(\cdot, t_n) - v^{[n]}\| \lesssim \tau^2$ established in \cite{BaoFS}, our estimates are improved in the sense that the leading error term as $\tau\to0^+$ is  now $\eps^2\tau^2$, which was numerically observed in \cite{BaoFS}.
The estimates in Theorem \ref{thm:semi} hold for higher order norms $\|\cdot\|_s$ ($s\leq m$) and the proof remains the same. The results are valid in higher dimensions $d=2,3$, independent of the aspect ratio of the  rectangular domain $\Omega$.
\end{remark}
\begin{remark} The second-order Strang splitting method is used to discretize the NKGE \eqref{eq:21} and it is straightforward to design the first-order Lie-Trotter splitting method \cite{Trotter} and fourth-order partitioned Runge-Kutta (PRK) splitting method \cite{BL, Geng}. Under appropriate assumptions of the exact solution, the improved uniform error bounds could be extended to the first-order Lie-Trotter splitting and the fourth-order PRK splitting method with improved uniform error bounds at $\varepsilon^2 \tau$ and $\varepsilon^2 \tau^{4}$, respectively.


\end{remark}

\subsection{Proof for Theorem \ref{thm:semi}}
The assumption (A) is equivalent to the regularity of $\psi(x, t)$ as $
\|\psi\|_{L^{\infty}\left([0, T_{\varepsilon}]; H^{m+1}\right)} \lesssim 1$.
Denote
\begin{equation}
F_t: \ \phi \mapsto e^{-it\langle\nabla\rangle}  F\big(e^{it\langle\nabla\rangle}\phi\big),\quad t \in \mathbb{R},
\label{eq:Ft}
\end{equation}
and we have the following estimates by the standard analysis for the local truncation error \cite{BCF,BaoFS}.
\begin{lemma}
\label{localeror}
For $0<\varepsilon\le 1$, the local error of the Strang splitting \eqref{S2} can be written as
\begin{equation}
\mathcal{E}^{n} :=  \mathcal{S}_{\tau}(\psi(t_n))- \psi(t_{n+1}) = \mathcal{F}(\psi(t_n))  + \mathcal{R}^n, \quad n = 0, 1, \cdots,
\label{eq:local}
\end{equation}
where
\begin{equation}\label{eq:mF-def}
\mathcal{F}(\psi(t_n)) 	=  \varepsilon^2 e^{i\tau\langle\nabla\rangle}\left(\tau F_{\tau/2}(\psi(t_n)) - \int^{\tau}_0 F_{\theta}(\psi(t_n))d \theta \right),
\end{equation}
and the following error bounds  hold under the assumption (A) with $m\ge1$,
\begin{equation}
\left\|\mathcal{F}(\psi(t_n))\right\|_1  \lesssim \varepsilon^2\tau^3, \quad \left\|\mathcal{R}^n\right\|_1\lesssim \varepsilon^4\tau^3.
\label{eq:localbound}
\end{equation}
\end{lemma}

 Under the  assumption (A), for $0<\tau\leq\tau_c$, we have the estimates \eqref{reg_semi} on the numerical solution $\psi^{[n]}$, which provide the control on the nonlinearity. Thus, we focus on the refined estimates in Theorem \ref{thm:semi}.
Introduce the numerical error function $e^{[n]}:=e^{[n]}(x)$ ($n=0,1,\cdots$) as
\begin{equation}
e^{[n]}:=\psi^{[n]}-\psi(t_n),
\end{equation}
and we have the error equation from  \eqref{S2} and \eqref{eq:local} as
\begin{align}
e^{[n+1]} & = \mathcal{S}_{\tau}(\psi^{[n]})	- \mathcal{S}_{\tau}(\psi(t_n)) + \mathcal{E}^{n}  = e^{i\tau\langle\nabla\rangle}e^{[n]} + W^n+ \mathcal{E}^{n},\quad n\ge0,
\label{eq:eg_semi}
\end{align}
where $W^n:=W^n(x)$ ($n=0,1,\cdots$) is given by
\begin{equation*}
W^n(x) = \varepsilon^2 \tau e^{i\frac{\tau}{2}\langle\nabla\rangle}\left(F\left(e^{i\frac{\tau}{2}\langle\nabla\rangle}\psi^{[n]}\right) -F\left(e^{i\frac{\tau}{2}\langle\nabla\rangle}\psi(t_n)\right) \right).
\end{equation*}
Under the assumption (A), we have from \eqref{eq:F} and the estimates on $\psi^{[n]}$ in \eqref{reg_semi}  that \begin{equation}
\left\|W^n(x)\right\|_1 \lesssim \varepsilon^2 \tau \left\|F\left(e^{i\frac{\tau}{2}\langle\nabla\rangle}\psi^{[n]}\right) -F\left(e^{i\frac{\tau}{2}\langle\nabla\rangle}\psi(t_n)\right)\right\|_1\lesssim\eps^2\tau\left\|e^{[n]}\right\|_1.
\label{eq:Wb}
\end{equation}
Based on \eqref{eq:eg_semi}, we obtain
\begin{equation}
e^{[n+1]} = e^{i(n+1)\tau\langle\nabla\rangle}e^{[0]} + \sum\limits_{k=0}^n e^{i(n-k)\tau\langle\nabla\rangle}\Big(W^k(x) + \mathcal{E}^k\Big),\quad 0\leq n\leq T_\eps/\tau-1.
\end{equation}
Noticing $e^{[0]} = 0$, \eqref{eq:local}, \eqref{eq:localbound} and  \eqref{eq:Wb}, we have the estimates for $0\leq n\leq T_\eps/\tau-1$,
\begin{equation}\label{eq:egrow}
\|e^{[n+1]}\|_1	\lesssim \varepsilon^2\tau^2 + \varepsilon^2 \tau  \sum_{k=0}^n \|e^{[k]}\|_1 + \| \sum\limits_{k=0}^n e^{i(n-k)\tau\langle\nabla\rangle}\mathcal{F}(\psi(t_k))\|_1.\end{equation}
Direct applications of \eqref{eq:localbound} and Gronwall's inequality lead to the uniform error estimates
$\|e^{[n+1]}\|_1\lesssim \tau^2$ ($0\leq n\leq T_\eps/\tau-1$) as shown in \cite{BFS}. To analyze the error more carefully, we shall employ the {\textbf{regularity compensation oscillation}} (RCO) technique \cite{BCF} to deal with the last term on the RHS of \eqref{eq:egrow}. The key idea is a summation-by-parts procedure combined with spectrum cut-off and phase cancellation.

The first step is a spectral projection on $\psi(t_k)$ such that only finite Fourier modes of $\psi(t_k)$ need to be considered and the projection error could be controlled by the regularity of $\psi(t_k)$. The second step is to apply the summation-by-parts formula for the low Fourier modes in a proper way, such that the phase can be cancelled for small $\tau$ (the terms of the type $\sum_{k=0}^n e^{i(n-k)\tau\langle\nabla\rangle}$) and an extra order of $\eps^2$ could be gained from the terms like $\mathcal{F}(\psi(t_k))-\mathcal{F}(\psi(t_{k+1}))$.

Now, we demonstrate our strategy in detail.  From the relativistic NLSE \eqref{eq:NLSE}, we find that
$\partial_t \psi(x, t) - i\langle\nabla\rangle \psi(x, t) = i\eps^2F(\psi(x,t))=O(\eps^2)$. Thus, in order to gain an extra order of $\eps^2$, instead of $\psi(x,t)$,  it is nature to consider the `twisted variable' given by
 \begin{equation}
 \phi(x, t) = e^{-it\langle\nabla\rangle}\psi(x, t), \quad t \geq 0,
 \label{eq:twist_def}
 \end{equation}
which satisfies the equation $\partial_t\phi(x,t)=\eps^2e^{-it\langle\nabla\rangle}F(e^{it\langle\nabla\rangle}\phi(x,t))$. Under the assumption (A), we have $\|\phi\|_{L^{\infty}([0, T_{\eps}]; H^{m+1})}\lesssim1$ and $\left\|\partial_t \phi\right\|_{L^{\infty}([0, T_{\eps}]; H^{m+1})} \lesssim \eps^2$ with
\begin{equation}
\|\phi(t_{n+1}) - \phi(t_{n})\|_{m+1} \lesssim \varepsilon^2 \tau, \quad 0 \leq n \leq T_\eps/\tau-1.
\label{eq:twist}
\end{equation}
The RCO technique will be used to force $\partial_t \phi(t)$ to appear with a gain of order $O(\eps^2)$ for the summation-by-parts procedure in $\sum_{k=0}^n e^{i(n-k)\tau\langle\nabla\rangle}\mathcal{F}(\psi(t_k))$. Then  small $\tau$ is required to control the accumulation of the frequency of the type $e^{i(n-k)\tau\langle\nabla\rangle}$.

{\bf Step 1}. As introduced in \cite{BCF}, we start with the choice of the cut-off parameter on the Fourier modes. Let $\tau_0 \in (0, 1)$ and choose $N_0 = 2\lceil 1/\tau_0\rceil \in \mathbb{Z}^+$ ($\lceil \cdot\rceil$ is the ceiling function) with $1/\tau_0 \leq N_0/2 < 1 + 1/\tau_0$.  Under the assumption (A), recalling $F_t$ \eqref{eq:Ft} and the operator $\langle\nabla\rangle^{-1}$, we have
\begin{equation}
\|F_t(e^{it_k\langle\nabla\rangle}\phi(t_k))\|_{m+2}\lesssim \|\phi(t_k)\|_{m+1}^3\lesssim 1,\quad t\in\mathbb{R},\quad 0\leq k\leq \frac{T_\eps}{\tau},
\end{equation}
and the following estimates hold by the standard Fourier projection properties for $s\in[0,m+1]$,
\begin{equation}
\|F_t(\phi(t_k))-P_NF_t(\phi(t_k))\|_{s}+\tau_0\|\phi(x,t_k)-P_{N_0}\phi(x,t_k)\|_s\lesssim \tau_0^{m+2-s}.
\end{equation}
Combing the above estimates, \eqref{eq:F}, \eqref{eq:mF-def} and assumption (A), we derive for $0\leq k\leq T_\eps/\tau$,
\begin{align}
&\|P_{N_0}\mathcal{F}(e^{i t_n\langle\nabla\rangle}(P_{N_0}\phi(t_k)))-
\mathcal{F}(e^{i t_k\langle\nabla\rangle}\phi(t_k))\|_1\nonumber\\
&\lesssim \eps^2\tau\tau_0^{m+1}+\eps^2\tau\|P_{N_0}\phi(t_k)-\phi(t_k)\|
\lesssim \eps^2\tau\tau_0^{m+1},
\end{align}
and \eqref{eq:egrow}  would imply for $0\leq n\leq T_\eps/\tau-1$,
\begin{equation}
\left\|e^{[n+1]}\right\|_1  \lesssim \tau_0^{m+1} +\varepsilon^2\tau^2+\varepsilon^2\tau\sum_{k=0}^n\left\|e^{[k]}\right\|_1 +\left\|\mathcal{L}^n\right\|_1, \label{eq:final2}
\end{equation}
where
\begin{equation}
 \mathcal{L}^n = \sum\limits_{k=0}^n e^{-i(k+1)\tau\langle\nabla\rangle}P_{N_0}\mathcal{F}(e^{i t_k\langle\nabla\rangle}(P_{N_0}\phi(t_k))).
\end{equation}

{\bf Step 2}. Now, we concentrate on the low Fourier modes term $\mathcal{L}^n$. Recalling the nonlinear function $ F(\cdot)$, we have the decomposition
\begin{equation}
F(\phi)=\sum_{q=1}^4F^q(\phi),\quad F^q(\phi)=i\langle\nabla\rangle^{-1}G^q(\phi),\quad q=1,2,3,4,
\end{equation}
with 
\begin{equation}
 G^1(\phi)=\frac{1}{8}\bar\phi^3,\; G^2(\phi)=\frac{3}{8}\phi\bar{\phi}^2,\;G^3(\phi)=\frac{3}{8}\phi^2\bar{\phi},\;G^4(\phi)=\frac{1}{8}{\phi}^3.
\end{equation}
For $\theta\in\mathbb{R}$ and $q=1,2,3,4$, introducing $F_\theta^q(\psi(t_k))=e^{-i\theta\langle\nabla\rangle}F^q(e^{i\theta\langle\nabla\rangle}\psi(t_k))$ and
\begin{equation}
\mathcal{F}^q(\psi(t_k))=\varepsilon^2 e^{i\tau\langle\nabla\rangle}\left(\tau F_{\tau/2}^q(\psi(t_n)) - \int^{\tau}_0 F^q_{\theta}(\psi(t_n))d \theta \right),
\end{equation}
recalling \eqref{eq:Ft} and \eqref{eq:mF-def}, we have
\begin{equation}\label{eq:Ln1}
\mathcal{L}^n=\sum_{q=1}^4\mathcal{L}_q^n,\quad \mathcal{L}_q^n = \sum\limits_{k=0}^n e^{-i(k+1)\tau\langle\nabla\rangle}P_{N_0}\mathcal{F}^q(e^{i t_k\langle\nabla\rangle}(P_{N_0}\phi(t_k))),\; 1\leq q\leq4.
\end{equation}
Since the estimates on $\mathcal{L}^n_q$ ($q=1,2,3,4$) are  the same, we only present the case for $\mathcal{L}^n_1$ ($0\leq n\leq T_\eps/\tau-1$).
For $l\in\mathcal{T}_{N_0}$, define the index set $\mathcal{I}_l^{N_0}$ associated to $l$ as
\begin{equation}
\mathcal{I}_l^{N_0}=\left\{(l_1,l_2,l_3)\ | \ l_1+l_2+l_3=l,\ l_1,l_2,l_3\in\mathcal{T}_{N_0}\right\},
\end{equation}
and  the following expansion  holds in view of $P_{N_0}\phi(t_k)=\sum_{l\in\mathcal{T}_{N_0}}\widehat{\phi}_l(t_k)e^{i\mu_l(x-a)}$,
\begin{align*}
& e^{-it_{k+1}\langle\nabla\rangle}P_{N_0}( e^{i\tau\langle\nabla\rangle}F_{\theta}^1(e^{i t_k\langle\nabla\rangle}P_{N_0}\phi(t_k))) \\
&=\quad \sum\limits_{l\in\mathcal{T}_{N_0}}\sum\limits_{(l_1,l_2,l_3)\in\mathcal{I}_l^{N_0}}\frac{i}{8\delta_l} \mathcal{G}^k_{l,l_1,l_2,l_3}(\theta) e^{i\mu_l(x-a)},
\end{align*}
where the coefficients  $\mathcal{G}^k_{l,l_1,l_2,l_3}(\theta)$ are functions of $\theta\in\mathbb{R}$ defined as
\begin{equation}
\mathcal{G}^k_{l,l_1,l_2,l_3}(\theta)
= e^{-i(t_k+\theta)\delta_{l,l_1,l_2,l_3}}\widehat{\phi}_{l_1}(t_k)\widehat{\phi}_{l_2}(t_k)\widehat{\phi}_{l_3}(t_k)
\label{eq:mGdef}
\end{equation}
with $\delta_{l,l_1,l_2,l_3} = \delta_l + \delta_{l_1} +\delta_{l_2} + \delta_{l_3}$ and $\delta_l = \sqrt{1+\mu_l^2}$ for $l \in \mathcal{T}_{N_0}$. Thus, we have
\begin{equation}\label{eq:remainder-dec}
\mathcal{L}_1^n  = \frac{i\eps^2}{8}\sum\limits_{k=0}^n
\sum\limits_{l\in\mathcal{T}_{N_0}}\sum\limits_{(l_1,l_2,l_3)\in\mathcal{I}_l^{N_0}}\frac{1}{\delta_l}\Lambda^k_{l,l_1,l_2,l_3}e^{i\mu_l(x-a)},
\end{equation}
where
\begin{align}
\Lambda^k_{l,l_1,l_2,l_3} & = - \tau \mathcal{G}^k_{l,l_1,l_2,l_3}(\tau/2) +\int_0^\tau\mathcal{G}^k_{l,l_1,l_2,l_3}(\theta)\,d \theta\nonumber\\
& = r_{l,l_1,l_2,l_3}e^{-it_k\delta_{l,l_1,l_2,l_3}}c^k_{l,l_1,l_2,l_3},
\label{eq:Lambdak}
\end{align}
with coefficients $c^k_{l,l_1,l_2,l_3}$ and $r_{l,l_1,l_2,l_3}$ given by
\begin{align}
c^k_{l,l_1,l_2,l_3}=& \ \widehat{\phi}_{l_1}(t_k) \widehat{\phi}_{l_2}(t_k)\widehat{\phi}_{l_3}(t_k),\label{eq:calFdef}\\
r_{l,l_1,l_2,l_3}= & \ -\tau e^{-i\tau\delta_{l,l_1,l_2,l_3}/2}+\int_0^\tau e^{-i\theta\delta_{l,l_1,l_2,l_3}}\,d \theta = O\left(\tau^3 (\delta_{l,l_1,l_2,l_3})^2\right).\label{eq:rest}
\end{align}
We only need consider the case $\delta_{l,l_1,l_2,l_3}\neq 0$ as $r_{l,l_1,l_2,l_3} = 0$ if $\delta_{l, l_1, l_2, l_3}= 0$. For $l\in\mathcal{T}_{N_0}$ and $(l_1,l_2,l_3)\in\mathcal{I}_l^{N_0}$, we have
\begin{equation}
|\delta_{l,l_1,l_2,l_3}|\leq 4\delta_{N_0/2}= 4 \sqrt{1+\mu_{N_0/2}^2} <  4\sqrt{1+\frac{4\pi^2(1+\tau_0)^2}{\tau_0^2(b-a)^2}},
\end{equation}
which implies when $0 < \tau \leq \alpha \frac{\pi (b-a)\tau_0}{2\sqrt{\tau_0^2(b-a)^2+4\pi^2(1+\tau_0)^2}}:=\tau_0^\alpha$ ($0 < \tau_0, \alpha <1$),  there holds 
\begin{equation}\label{eq:cfl}
\frac{\tau}{2}|\delta_{l,l_1,l_2,l_3}| \leq \alpha\pi.
\end{equation}
Denoting $S^n_{l,l_1,l_2,l_3}=\sum_{k=0}^n e^{-it_k\delta_{l,l_1,l_2,l_3}}$ ($n\ge0$), for  $0< \tau \leq \tau_0^\alpha$, we then obtain from \eqref{eq:cfl} that
\begin{equation}\label{eq:Sbd}
|S^n_{l,l_1,l_2,l_3}|\leq \frac{1}{|\sin(\tau \delta_{l,l_1,l_2,l_3}/2)|}\leq\frac{C}{\tau|\delta_{l,l_1,l_2,l_3}|},\quad C = \frac{2\alpha\pi}{\sin(\alpha \pi)}\quad \forall n\ge0.
\end{equation}
Using summation by parts,
we find from \eqref{eq:Lambdak} that
\begin{equation}
\sum_{k=0}^n\Lambda^k_{l,l_1,l_2,l_3}=r_{l,l_1,l_2,l_3}\big[\sum_{k=0}^{n-1}S^k_{l,l_1,l_2,l_3} (c^k_{l,l_1,l_2,l_3}-c^{k+1}_{l,l_1,l_2,l_3})+S^n_{l,l_1,l_2,l_3}c^n_{l,l_1,l_2,l_3}\big],
\label{eq:lambdasum}
\end{equation}
with
\begin{align}
&c^k_{l,l_1,l_2,l_3}-c^{k+1}_{l,l_1,l_2,l_3}\nn\\
& = (\widehat{\phi}_{l_1}(t_k)-\widehat{\phi}_{l_1}(t_{k+1})) \widehat{\phi}_{l_2}(t_k)\widehat{\phi}_{l_3}(t_k) + \widehat{\phi}_{l_1}(t_{k+1})(\widehat{\phi}_{l_2}(t_k)-\widehat{\phi}_{l_2}(t_{k+1}))\widehat{\phi}_{l_3}(t_k)\nn \\
&\;\;\;\;\; + \widehat{\phi}_{l_1}(t_{k+1}) \widehat{\phi}_{l_2}(t_{k+1})(\widehat{\phi}_{l_3}(t_k)-\widehat{\phi}_{l_3}(t_{k+1})).\label{eq:cksum}
\end{align}

Combining \eqref{eq:rest}, \eqref{eq:Sbd}, \eqref{eq:lambdasum} and \eqref{eq:cksum}, we have
\begin{align}
\left|\sum_{k=0}^n\Lambda^k_{l,l_1,l_2,l_3}\right|
\lesssim & \  \tau^2|\delta_{l,l_1,l_2,l_3}|\sum\limits_{k=0}^{n-1}\bigg(
\left|\widehat{\phi}_{l_1}(t_k)-\widehat{\phi}_{l_1}(t_{k+1})\right|\left|\widehat{\phi}_{l_2}(t_k)\right| \left|\widehat{\phi}_{l_3}(t_k)\right|\nn\\
&\ +\left|\widehat{\phi}_{l_1}(t_{k+1})\right|\left|\widehat{\phi}_{l_2}(t_k)-\widehat{\phi}_{l_2}(t_{k+1})\right| \left|\widehat{\phi}_{l_3}(t_k)\right|\nonumber\\
&\ +\left|\widehat{\phi}_{l_1}(t_{k+1})\right|\left|\widehat{\phi}_{l_2}(t_{k+1})\right| \left|\widehat{\phi}_{l_3}(t_k)-\widehat{\phi}_{l_3}(t_{k+1})\right|\bigg)\nn\\
&\ + \tau^2|\delta_{l,l_1,l_2,l_3}| \left|\widehat{\phi}_{l_1}(t_n)\right|\left|\widehat{\phi}_{l_2}(t_n)\right| \left|\widehat{\phi}_{l_3}(t_n)\right|.\label{eq:sumlambda}
\end{align}
For $l\in\mathcal{T}_{N_0}$ and $(l_1,l_2,l_3)\in\mathcal{I}_l^{N_0}$, there holds
\begin{equation}
|\delta_{l,l_1,l_2,l_3}| \leq \big(1+(\sum^3_{j=1}\mu_{l_j})^2\big)^{1/2} + \sum^3_{j=1} \sqrt{1+\mu_{l_j}^2} \lesssim \prod_{j=1}^3\sqrt{1+\mu_{l_j}^2},
\label{eq:mlbd}
\end{equation}
Based on \eqref{eq:remainder-dec}, \eqref{eq:sumlambda} and \eqref{eq:mlbd}, noticing $\delta_l = \sqrt{1+\mu_l^2}$, we have
\begin{align}
&\left\|\mathcal{L}_1^n\right\|^2_1  \nn \\
& = \ \varepsilon^4
\sum\limits_{l\in\mathcal{T}_{N_0}}\left|\sum\limits_{(l_1,l_2,l_3)\in\mathcal{I}_l^{N_0}}\sum\limits_{k=0}^n\Lambda^k_{l,l_1,l_2,l_3}\right|^2\nn \\
& \lesssim \ \varepsilon^4\tau^4
\bigg\{\sum_{l\in\mathcal{T}_{N_0}}\bigg(\sum\limits_{(l_1,l_2,l_3)\in\mathcal{I}_l^{N_0}}\left|\widehat{\phi}_{l_1}(t_n)\right|\left|\widehat{\phi}_{l_2}(t_n)\right| \left|\widehat{\phi}_{l_3}(t_n)\right|\prod_{j=1}^3\sqrt{1+\mu_{l_j}^2}\bigg)^2\nn \\
&\;\;\;\; + n \sum\limits_{k=0}^{n-1}\sum_{l\in\mathcal{T}_{N_0}}\bigg[\bigg(\sum\limits_{(l_1,l_2,l_3)\in\mathcal{I}_l^{N_0}}
\left|\widehat{\phi}_{l_1}(t_k)-\widehat{\phi}_{l_1}(t_{k+1})\right|\left|\widehat{\phi}_{l_2}(t_k)\right| \left|\widehat{\phi}_{l_3}(t_k)\right|\prod_{j=1}^3\sqrt{1+\mu_{l_j}^2}\bigg)^2\nn\\
& \;\;\;\; +\bigg(\sum\limits_{(l_1,l_2,l_3)\in\mathcal{I}_l^{N_0}}
\left|\widehat{\phi}_{l_1}(t_{k+1})\right|\left|\widehat{\phi}_{l_2}(t_k)-\widehat{\phi}_{l_2}(t_{k+1})\right| \left|\widehat{\phi}_{l_3}(t_k)\right|\prod_{j=1}^3\sqrt{1+\mu_{l_j}^2}\bigg)^2\nn\\
& \;\;\;\; +\bigg(\sum\limits_{(l_1,l_2,l_3)\in\mathcal{I}_l^{N_0}}
\left|\widehat{\phi}_{l_1}(t_{k+1})\right|\left|\widehat{\phi}_{l_2}(t_{k+1})\right| \left|\widehat{\phi}_{l_3}(t_k)-\widehat{\phi}_{l_3}(t_{k+1})\right|\prod_{j=1}^3\sqrt{1+\mu_{l_j}^2}\bigg)^2\bigg]\bigg\}.
\label{eq:sumlambda-2}
\end{align}
In order to estimate the sum on the RHS of above inequality, e.g. for the first term on the RHS, we  use the auxiliary function $\xi(x)=\sum_{l\in\mathbb{Z}}\sqrt{1+\mu_{l_j}^2}\left|\widehat{\phi}_l(t_n)\right|e^{i\mu_l(x-a)}$, where $\xi(x)\in H^{m}(\Omega)$ implied by assumption (A) and $\|\xi\|_{H^s}\lesssim\|\phi(t_n)\|_{H^{s+1}}$ ($s \leq m$). Expanding $|\xi(x)|^2\xi(x)=\sum\limits_{l\in\mathbb{Z}}\sum\limits_{l_1-l_2+l_3=l, l_j\in\mathbb{Z}} \prod_{j=1}^3\left(\sqrt{1+\mu_{l_j}^2}\left|\widehat{\phi}_{l_j}(t_n)\right|\right)e^{i\mu_l(x-a)}$, we could obtain
\begin{align}
&\sum_{l\in\mathcal{T}_{N_0}} \bigg(\sum\limits_{(l_1,l_2,l_3)\in\mathcal{I}_l^{N_0}}\left|\widehat{\phi}_{l_1}(t_n)\right|\left|\widehat{\phi}_{l_2}(t_n)\right| \left|\widehat{\phi}_{l_3}(t_n)\right|\prod_{j=1}^3\sqrt{1+\mu_{l_j}^2}\bigg)^2\nn \\
& \  \leq \left\||\xi(x)|^2\xi(x)\right\|^2 \lesssim \left\|\xi(x)\right\|_1^6\lesssim \left\|\phi(t_n)\right\|_2^6 \lesssim1.
\end{align}
Thus, in light of \eqref{eq:twist}, we could estimate each term in \eqref{eq:sumlambda-2} similarly as
\begin{align}
\|\mathcal{L}_1^n\|_1& \ \lesssim \varepsilon^4\tau^4 \bigg[\left\|\phi(t_n)\right\|_2^6+n\sum\limits_{k=0}^{n-1}
\left\|\phi(t_k) - \phi(t_{k+1})\right\|_2^2(\left\|\phi(t_k)\right\|_2 + \left\|\phi(t_{k+1})\right\|_2)^4\bigg]\nn \\
&\ \lesssim  \varepsilon^4\tau^4+n^2\varepsilon^4\tau^4 (\varepsilon^2\tau)^2
\lesssim \varepsilon^4\tau^4,\quad 0\leq n\leq T_\eps/\tau-1.
\label{eq:est-l2}
\end{align}
The same estimates could be established for $\mathcal{L}_q^n$ ($q=2,3,4$)
and \eqref{eq:final2}  together with \eqref{eq:Ln1}  implies
\begin{equation}
\left\|e^{[n+1]}\right\|_1 \lesssim \tau_0^{m+1} +\eps^2\tau^2+\eps^2\tau\sum_{k=0}^n\left\|e^{[k]}\right\|_1,\quad 0\leq n \leq T_\eps/\tau-1.	
\end{equation}
Discrete Gronwall's inequality yields
\begin{equation}
\left\|e^{[n+1]}\right\|_1 \lesssim \eps^2\tau^2 + \tau_0^{m+1},\quad 0\leq n \leq T_\eps/\tau-1,
\end{equation}
and  the error bound \eqref{eq:error_semi} follows in  view of \eqref{uxtpsit} and \eqref{uxtpsita}.
\begin{remark} Similar results in Theorem \ref{thm:semi} have been previously obtained for the time-splitting method applied to the long-time dynamics of nonlinear Schr\"odinger equation with weak nonlinearity \cite{CMTZ}, where the periodicity of the free Schr\"odinger  operator plays an important role and the time step size has to be an integer fraction of the period. Thus, the results and analysis in \cite{CMTZ} are difficult to extend to the higher dimensional rectangular domain with irrational aspect ratio, and/or the general time step sizes. The presented RCO based approach does not depend on the periodicity of the free relativistic Schr\"odinger operator. It is easy to check our analysis works for the higher dimensional cases and allows general time step sizes.
\end{remark}

\section{Full-discretization and improved uniform error bounds}
In this section, we present the practical full-discretization for the NKGE \eqref{eq:21} by the Fourier pseudospectral method in space and establish the improved uniform error bounds.

\subsection{Full-discretization by Fourier pseudospectral method}
Let $N$ be an even  positive integer and define the spatial mesh size $h = (b - a)/N$, then the grid points are chosen as
\begin{equation}
x_j := a + jh,\quad j \in \mathcal{T}^0_N=\{j~|~j = 0, 1, \ldots, N\}.
\end{equation}

Let $\psi_j^n$ be the numerical approximation of $\psi(x_j,t_n)$
for $j\in \mathcal{T}^0_N$ and $n\ge0$ and denote $\psi^n=(\psi_0^n, \psi_1^n,\ldots, \psi_N^n)^T\in \mathbb{C}^{N+1}$ for $n=0,1,\ldots$. Then a time-splitting
Fourier pseudospectral (TSFP) method for discretizing the relativistic NLSE  \eqref{eq:NLSE} via \eqref{S2} with a Fourier pseudospectral discretization in space is given as
\begin{equation}
\label{psifull}
\begin{split}
&\psi^{(1)}_j=\sum_{l \in \mathcal{T}_N} e^{i\frac{\tau\delta_l}{2}}\;\widetilde{(\psi^n)}_l\; e^{i\mu_l(x_j-a)},  \\
&\psi^{(2)}_j=\psi^{(1)}_j+ \eps^2\tau\, F_j^n, \qquad F_j^n=i\sum_{l \in \mathcal{T}_N} \frac{1}{\delta_l} \widetilde{\left(G(\psi^{(1)})\right)}_l\; e^{i\mu_l(x_j-a)}, \\
&\psi^{n+1}_j=\sum_{l \in \mathcal{T}_N} e^{i\frac{\tau\delta_l}{2}} \; \widetilde{\left(\psi^{(2)}\right)}_l\; e^{i\mu_l(x_j-a)},\quad j \in \mathcal{T}^0_N, \quad n=0,1,\ldots,
\end{split}
\end{equation}
where $\delta_l = \sqrt{1+\mu_l^2}$ for $l \in \mathcal{T}_N$, $\psi^{(k)}=(\psi_0^{(k)}, \psi_1^{(k)},\ldots$, $\psi_N^{(k)})^T\in \mathbb{C}^{N+1}$ for
$k=1$, $2$, $G(\psi^{(1)}):=(G(\psi^{(1)}_0), G(\psi^{(1)}_2), \ldots, G(\psi^{(1)}_N))^T\in \mathbb{R}^{N+1}$ and
\[\psi_j^0 = u_0(x_j)-i\sum_{l \in \mathcal{T}_N}\frac{\widetilde{(u_1)}_l}{\delta_l} e^{i\mu_l(x_j-a)}, \quad j \in \mathcal{T}^0_N. \]

Let $u^n_j$ and $v^n_j$ be the approximations of $u(x_j, t_n)$ and $v(x_j, t_n)$, respectively, for $j\in \mathcal{T}^0_N$ and $n\ge0$, and denote $u^n=(u_0^n, u_1^n,\ldots, u_N^n)^T\in \mathbb{R}^{N+1}$ and  $v^n=(v_0^n, v_1^n,\ldots, v_N^n)^T\in \mathbb{R}^{N+1}$  for $n=0,1,\ldots$. Combining \eqref{uxtpsita} and \eqref{psifull}, we could obtain the
full-discretization of the NKGE \eqref{eq:21} by the TSFP method as
\begin{equation}
\label{ufull}
\begin{split}
&u_j^{n+1}=\frac{1}{2}\left(\psi_j^{n+1}+\overline{\psi_j^{n+1}}\right),\\
&v_j^{n+1} = \frac{i}{2}\sum_{l \in \mathcal{T}_N}\delta_l\big[\widetilde{(\psi^{n+1})}_l-
\widetilde{(\overline{\psi^{n+1}})}_l\big]\;
e^{i\mu_l (x_j-a)},
\end{split}
\qquad j \in \mathcal{T}^0_N, \quad n\ge 0,
\end{equation}
with $u_j^0=u_0(x_j)$ and  $v_j^0 = u_1(x_j)$ for $j \in \mathcal{T}^0_N$.

\subsection{Improved uniform error bounds}Let $u^n$ and $v^n$ be the numerical approximations obtained from the TSFP \eqref{psifull}--\eqref{ufull}. From the analysis in \cite{BaoFS}, under the assumption (A),  for $0<\tau\leq\tau_c,0< h\leq h_c$ ($\tau_c$, $h_c$ are constants independent of $\eps$), there exists a constant $M>0$ depending on $T$, $\|u_0\|_{m+1}$, $\|u_1\|_m$, $\|u\|_{L^{\infty}([0, T_\eps]; H^m)}$ and $\|\partial_t u\|_{L^{\infty}([0, T_\eps]; H^m)}$ such that the numerical solution satisfies
\begin{equation}\label{reg_full}
\left\|I_N u^n\right\|_{m+1}^2 + \left\|I_N v^n\right\|_{m}^2 \leq M,\; \text{or equivalently}\;
\|I_N\psi^n\|_{m+1}^2\leq M \quad 0 \leq n \leq \frac{T_\eps}{\tau}.
\end{equation}
Then we have the improved uniform error bounds for the full-discretization.

\begin{theorem}
\label{thm:full}Under the assumption (A), there exist $h_0 > 0$ and $0 < \tau_0 < 1$ sufficiently small and independent of $\varepsilon$ such that, for any $0 < \varepsilon \leq 1$,  when $0 < h \leq h_0$ and $0< \tau < \alpha \frac{\pi (b-a)\tau_0}{2\sqrt{\tau^2_0(b-a)^2 + 4\pi^2(1+\tau_0^2)}}$ for a fixed constant $\alpha \in (0, 1)$, we have the following improved uniform error estimates
\begin{equation}
\|u(\cdot, t_n) - I_Nu^n\|_1 + \|\partial_t u(\cdot, t_n) - I_Nv^n\| \lesssim h^{m} + \varepsilon^2\tau^2 + \tau^{m+1}_0,\ 0 \leq n \leq \frac{T/\varepsilon^2}{\tau}.
\label{eq:error_full}
\end{equation}
In particular, if the exact solution is sufficiently smooth, e.g. $u, \partial_t u \in H^{\infty}$, the improved uniform error bounds for sufficiently small $\tau$ could be
\begin{equation}
\|u(\cdot, t_n) - I_Nu^n\|_1 + \|\partial_t u(\cdot, t_n) - I_Nv^n\| \lesssim h^{m} + \varepsilon^2\tau^2,\ 0 \leq n \leq \frac{T/\varepsilon^2}{\tau}.
\end{equation}
\end{theorem}

\emph{Proof.} It suffices to consider the numerical approximation $\psi^n$ to the solution of the relativistic NLSE  \eqref{eq:NLSE}. Recalling the semi-discrete-in-time approximation $\psi^{[n]}$ ($0 \leq n \leq \frac{T/\eps^2}{\tau}$) given by the scheme \eqref{S2}-\eqref{uxtpsita}, under the assumptions of Theorem \ref{thm:full}, we have the estimates in Theorem \ref{thm:semi}, \eqref{reg_semi} and \eqref{reg_full}, which directly yield
\begin{equation}
\label{psip}
\left\|\psi^{[n]}-P_N\psi^{[n]}\right\|_1 \lesssim  h^{m},\quad
\left\|\psi(\cdot, t_n) - \psi^{[n]}\right\|_1 \lesssim \varepsilon^2\tau^2+\tau_0^{m+1},\quad 0\leq n\leq \frac{T_\eps}{\tau}.
\end{equation}
Since $\psi(\cdot, t_n)-I_N\psi^n = \psi(\cdot, t_n)-\psi^{[n]}+\psi^{[n]}-P_N\psi^{[n]} + P_N\psi^{[n]} - I_N\psi^n$, we derive that
\begin{equation}
\label{sp}
\|\psi(\cdot, t_n)-I_N\psi^n\|_1\leq \|P_N\psi^{[n]} - I_N\psi^n\|_1+C_1 (\varepsilon^2\tau^2+\tau_0^{m+1}+h^m).
\end{equation}
As a result, it remains to establish the estimates  on the error function $e^n:=e^n(x)\in Y_N$  given as
\[
e^n:=P_N\psi^{[n]} - I_N\psi^n,\quad 0 \leq n \leq \frac{T/\eps^2}{\tau}.\]
From  \eqref{S2} and \eqref{psifull}, we get
\begin{align*}
&I_N\psi^{n+1} = e^{i\tau\langle\nabla\rangle}I_N\psi^{n}+i\eps^2 \tau \langle \nabla \rangle^{-1}e^{i\tau\langle\nabla\rangle/2}I_N(G(e^{i\tau\langle\nabla\rangle/2}I_N\psi^{n})),\\
&P_N\psi^{[n+1]} = e^{i\tau\langle\nabla\rangle}P_N\psi^{[n]}+i\eps^2 \tau \langle \nabla \rangle^{-1}e^{i\tau\langle\nabla\rangle/2}P_N(G(e^{i\tau\langle\nabla\rangle/2}\psi^{[n]})),
\end{align*}
which lead to
\begin{align}\label{eq:errg:f}
e^{n+1}=e^{i\tau\langle\nabla\rangle}e^n+i\eps^2 \tau \langle \nabla \rangle^{-1}e^{i\tau\langle\nabla\rangle/2}\left(P_NG(\psi^{\langle1\rangle})-I_NG(\psi^{(1)})\right),
\end{align}
with $\psi^{\langle1\rangle}=e^{i\tau\langle\nabla\rangle/2}\psi^{[n]}$ and $\psi^{(1)}=e^{i\tau\langle\nabla\rangle/2}I_N\psi^{n}$.
Hence, combining the bounds \eqref{reg_semi} and \eqref{reg_full}, we have $\|G(\psi^{\langle1\rangle})\|_{m+1}+\|G(\psi^{(1)})\|_{m+1}\lesssim 1$ and
\begin{equation}
\|G(\psi^{\langle1\rangle})-G(\psi^{(1)})\|\lesssim \|\psi^{\langle1\rangle}-\psi^{(1)}\|\lesssim \|\psi^{[n]}-I_N\psi^n\|\lesssim h^{m+1}+\|e^n\|.
\end{equation}
To summarize, noticing $\|P_NG(\psi^{\langle1\rangle})-I_NG(\psi^{(1)})\|\leq \|P_N(G(\psi^{(1)})) - I_N(G(\psi^{(1)}))\|+\|P_N(G(\psi^{\langle1\rangle}))-P_N(G(\psi^{(1)}))\|\lesssim h^{m+1}+ \|G(\psi^{\langle1\rangle})-G(\psi^{(1)})\|$, we  could obtain from \eqref{eq:errg:f} that
\begin{align*}
\|e^{n+1}\|_1 \leq& \|e^n\|_1+\eps^2\tau\|P_NG(\psi^{\langle1\rangle})-I_NG(\psi^{(1)})\| \\
\leq& \|e^n\|_1 + C_1\eps^2\tau h^{m+1}+C_2\eps^2\tau\|e^n\|,\quad 0\leq n\leq T_\eps/\tau-1,
\end{align*}
where $C_1,C_2$ are constants independent of $\eps,h,\tau,n,\tau_0$. Since
$e^0=P_Nu_0-I_Nu_0-i\langle\nabla\rangle^{-1}(P_Nu_1-I_Nu_1)$, we have
$\|e^0\|_1\lesssim h^m$ and  discrete Gronwall's inequality  implies $\|e^{n+1}\|\lesssim h^{m+1}$ ($0\leq n\leq T_\eps/\tau-1$). Combining the above estimtates with \eqref{sp}, we derive
\[\left\|\psi(\cdot, t_n) - I_N\psi^n\right\|_1 \lesssim h^{m} + \varepsilon^2\tau^2+\tau^{m+1}_0,\quad 0\leq n\leq T_\eps/\tau.\]
Recalling \eqref{ufull}, we obtain error bounds for $u^n$ and $v^n$  ($0\leq n\leq T_\eps/\tau$) as
\begin{align*}
\|u(\cdot, t_n)-I_Nu^n\|_1 & = \frac{1}{2}\left\|\psi(\cdot, t_n)+\overline{\psi(\cdot, t_n)}-
I_N\psi^n - I_N \overline{\psi^n}\right\|_1\\
&\le \|\psi(\cdot, t_n)-I_N\psi^n\|_1 \lesssim h^{m} + \varepsilon^2\tau^2 + \tau_0^{m+1},\\
 \|v(\cdot, t_n)-I_Nv^n\| &=\frac{1}{2}\|\langle \nabla\rangle(\psi(\cdot, t_n)-\overline{\psi(\cdot, t_n)})-\langle \nabla\rangle (I_N\psi^n - I_N\overline{\psi^n})\| \\
&\le \|\psi(\cdot, t_n)-I_N\psi^n\|_1 \lesssim h^{m} + \varepsilon^2\tau^2 + \tau_0^{m+1},
\end{align*}
which show \eqref{eq:error_full} and the proof for Theorem \ref{thm:full} is completed.
\hfill $\square$ \bigskip
\begin{remark}Through the proof of Theorem \ref{thm:full}, it is not difficult to see the spatial error estimates of $u(\cdot, t_n) - I_Nu^n$ in  $L^2$ norm can be improved to $h^{m+1}$.
\end{remark}

\section{Extensions}
In this section, we discuss the extensions of the time-splitting method and corresponding error estimates to the complex NKGE with a general power nonlinearity and an oscillatory complex NKGE which propagates waves with wavelength at $O(\eps^{2p})$ in time.
\subsection{To the complex NKGE with a general power nonlinearity}
Consider the following complex NKGE with a general power nonlinearity
\begin{equation}
\label{eq:GWNE}
\left\{
\begin{aligned}
&\partial_{tt}u(\bx, t)-\Delta u({\bx}, t)+ u(\bx, t)+\varepsilon^{2p} |u({\bx}, t)|^{2p} u({\bx}, t)=0,\quad \bx \in \Omega,\quad t > 0,\\
&u(\bx, 0) = u_0(\bx), \quad \partial_t u(\bx, 0) = u_1(\bx),\quad{\bx} \in \Omega.
\end{aligned}\right.
\end{equation}
Here, $u := u(\bx, t)$ is a complex-valued scalar field, $p \in \mathbb{N}^+$ is the power index, and the initial data $u_0(\bx)$ and $u_1(\bx)$ are two given complex-valued functions which are independent of $\eps$. The domain $\Omega$ and periodic boundary conditions are given the same as those in \eqref{eq:WNE}. The local/global well-posedness and scattering properties of the Cauchy problem \eqref{eq:GWNE} have been widely studied in the literature and references therein \cite{GV,KS,MS,MK,OTT,Pecher,Tao}. From the analytical results, the life-span of a smooth solution to the complex NKGE \eqref{eq:GWNE} is at least $O(\eps^{-2p})$.

For simplicity of notations, we only show the numerical scheme in 1D under the periodic boundary condtions. Similarly, introducing $v(x, t) = \partial_t u(x, t)$ and
\begin{equation}
\label{etapm}
\eta_{\pm}(x, t) =  u(x, t) \mp i\ \langle \nabla \rangle^{-1}v(x, t), \qquad a\le x\le b, \quad t\ge0,
\end{equation}
and denoting $f(\varphi)=|\varphi|^{2p}\varphi$, then the complex NKGE \eqref{eq:GWNE} can be reformulated into the following coupled relativistic NLSEs:
\begin{equation}
\label{eq:NLS_com}
\left\{
\begin{aligned}
&i\partial_t \eta_{\pm}\pm\langle\nabla \rangle \eta_\pm \pm \varepsilon^{2p} \langle\nabla \rangle^{-1}f\left(\frac{1}{2}\eta_++\frac{1}{2}\eta_-\right) = 0,\\
&\eta_\pm(t=0)=u_0\mp i\ \langle \nabla \rangle^{-1}v_0.
\end{aligned}\right.
\end{equation}

Let $\eta^n_{\pm,j}$ be the approximations of $\eta_{\pm}(x_j, t_n)$ for $j \in \mathcal{T}^0_N$ and $n \geq 0$, and denote $\eta_{\pm}^n=(\eta^n_{\pm,0},\eta^n_{\pm,1},\ldots,\eta^n_{\pm,N})^T\in \mathbb{C}^{n+1}$ as the solution at $t_n = n \tau$. Similar to the NKGE with cubic nonlinearity, the second-order time-splitting Fourier pseudospectral (TSFP) discretization for the relativistic NLSE \eqref{eq:NLS_com} is given by
\begin{equation}
\begin{split}
&\eta^{(1)}_{\pm,j}=\sum_{l \in \mathcal{T}_N} e^{\pm~i\frac{\tau\delta_l}{2}}\;\widetilde{(\eta_{\pm}^n)}_l \; e^{i\mu_l(x_j-a)},  \\
&\eta^{(2)}_{\pm,j}=\eta_{\pm,j}^{(1)}\pm\varepsilon^{2p}\tau f_j^n, \\
&\eta^{n+1}_{\pm,j}=\sum_{l \in \mathcal{T}_N} e^{\pm~i\frac{\tau\delta_l}{2}}\; \widetilde{(\eta_{\pm}^{(2)})}_l\; e^{i\mu_l(x_j-a)},
\end{split}\qquad j \in \mathcal{T}^0_N,\quad  n\ge0,
\label{eq:gpsi}
\end{equation}
with
\begin{equation*}
\eta_{\pm,j}^0=u_0(x_j)\mp i\sum_{l\in\mathcal{T}_N}\frac{1}{\delta_l}\widetilde{(v_0)}_le^{i\mu_l(x_j-a)}, \quad f_j^n=i\sum_{l \in \mathcal{T}_N} \frac{1}{\delta_l} \; \widetilde{\left(f((\eta_{+}^n+\eta_{-}^n)/2)\right)}_l \; e^{i\mu_l(x_j-a)}.\end{equation*}
Then $u^{n+1}_j$ and $v^{n+1}_j$ ($j \in \mathcal{T}^0_N, n\ge0$) which are approximations of $u(x_j, t_{n+1})$ and $v(x_j, t_{n+1})$, respectively, can be recovered by
\begin{equation}
\label{eq:gufull}
u_j^{n+1}=\frac{1}{2}\left(\eta_{+,j}^{n+1}+\eta_{-,j}^{n+1}\right),\;v_j^{n+1}=\frac{i}{2}\sum_{l \in \mathcal{T}_N}\delta_l\Big(\widetilde{(\eta_{+}^{n+1})}_l-\widetilde{(\eta_{-}^{n+1})}_l\Big)
e^{i\mu_l (x_j-a)}.
\end{equation}

We assume the exact solution $u:=u(x, t)$ of the NKGE \eqref{eq:GWNE} up to the time  $T_{\varepsilon, p} = T/\varepsilon^{2p}$ ($T>0$ fixed):
\begin{equation*}
{\rm(B)}\qquad
\|u\|_{L^{\infty}\left([0, T_{\eps, p}]; H^{m+1}\right)} \lesssim 1,\quad \;\; \|\partial_t u\|_{L^{\infty}\left([0, T_{\eps, p}]; H^{m}\right)} \lesssim 1,\quad m\ge1,
\end{equation*}
 then the following improved uniform error bounds for the TSFP method \eqref{eq:gpsi}--\eqref{eq:gufull} could be established up to the  time $T_{\eps,p}$.
\begin{theorem}\label{thm:general}
Let $u^n$ and $v^n$ be the numerical approximations obtained from the TSFP \eqref{eq:gpsi}--\eqref{eq:gufull}. Under the assumption (B), there exist $h_0 > 0$ and $0 < \tau_0 < 1$ sufficiently small and independent of $\varepsilon$ such that, for any $0<\eps\leq1$, when $0 < h \leq h_0$ and $0< \tau \leq\alpha \tau_0$ for some fixed constant $\alpha>0$, we have the following improved uniform error estimates
\begin{equation}
\|u(\cdot, t_n) - I_Nu^n\|_1 + \|\partial_t u(\cdot, t_n) - I_Nv^n\| \lesssim h^{m} + \eps^{2p}\tau^2 + \tau^{m+1}_0,\ 0 \leq n \leq \frac{T/\varepsilon^{2p}}{\tau}.
\end{equation}
In particular, if the exact solution is sufficiently smooth, e.g. $u, \partial_t u \in H^{\infty}$, the uniform improved error bounds for sufficiently small $\tau$ could be
\begin{equation}
\|u(\cdot, t_n) - I_Nu^n\|_1 + \|\partial_t u(\cdot, t_n) - I_Nv^n\| \lesssim h^{m} + \eps^{2p}\tau^2,\ 0 \leq n \leq \frac{T/\varepsilon^{2p}}{\tau}.
\end{equation}
\end{theorem}

\subsection{To an oscillatory complex NKGE}

Introducing a re-scale in time
\begin{equation}
t = \frac{r}{\eps^{2p}} \Leftrightarrow	 r = \eps^{2p} t, \quad \nu(\bx, r) = u(\bx, t),
\end{equation}
 the NKGE \eqref{eq:GWNE} could be reformulated into the following oscillatory complex NKGE
\begin{equation}
\label{eq:HOE}
\left\{
\begin{aligned}
& \eps^{2p} \partial_{rr}\nu(\bx, r) - \frac{1}{\eps^{2p}}\Delta \nu(\bx, r)+ \frac{1}{\eps^{2p}}\nu(\bx, r)+|\nu(\bx, r)|^{2p}\nu(\bx, r) = 0,\ \bx \in \Omega,\ r > 0,\\
&\nu(\bx, 0) = u_0({\bx}),\quad\partial_r \nu(\bx,0) = \frac{1}{\eps^{2p}}u_1({\bx}),\quad{\bx} \in \Omega.
\end{aligned}\right.
\end{equation}
The solution of the oscillatory NKGE \eqref{eq:HOE} propagates waves with amplitude at $O(1)$, wavelength at $O(1)$ and $O(\eps^{2p})$ in space and time, respectively, and wave velocity at $O(\eps^{-{2p}})$. Denote $\mu(\bx, r) = \partial_r \nu(\bx, r)$, by taking the time step $\kappa = \eps^{2p}\tau$, then the improved error bounds on the time-splitting methods  (see Remark \ref{remark:split}) for the long-time problem could be extended to the oscillatory complex NKGE \eqref{eq:HOE} up to the fixed time $T$.
\begin{theorem}
\label{thm:HOE}
Let $\nu^n$ and $\mu^n$ be the numerical approximations obtained from the TSFP method. Assume the exact solution $\nu$ of the oscillatory complex NKGE \eqref{eq:HOE} satisfies for some $m\ge1$:
\begin{equation*}
\begin{split}
&\nu \in  \  L^\infty\left([0, T]; H^{m+1}\right), \quad
\partial_r \nu\in L^\infty\left([0, T]; H^{m}\right),\\
&\|\nu\|_{L^{\infty}\left([0, T]; H^{m+1}\right)} \lesssim 1,\quad \;\; \|\partial_r \nu\|_{L^{\infty}\left([0, T]; H^{m}\right)} \lesssim \frac{1}{\eps^{2p}},
\end{split}
\end{equation*}
 there exist $h_0 > 0$ and $0 < \kappa_0 < 1$ sufficiently small and independent of $\varepsilon$ such that, for any $0 < \varepsilon \leq 1$, when the mesh size $0 < h \leq h_0$ and the time step $0< \kappa \leq\alpha \kappa_0\eps^{2p}$ for some fixed constant $\alpha>0$, we have the following improved error estimates
\begin{equation}
\|\nu(\cdot, r_n) - I_N\nu^n\|_1 + \eps^{2p}\|\partial_r \nu(\cdot, r_n) - I_N\mu^n\| \lesssim h^{m} +  \frac{\kappa^2}{\eps^{2p}} + \kappa_0^{m+1},\ 0 \leq n \leq \frac{T}{\kappa}.
\label{eq:HOE_error}
\end{equation}
In particular, if the exact solution is sufficiently smooth, e.g. $\nu, \partial_r \nu \in H^{\infty}$,  the improved error bounds for sufficiently small $\kappa$ could be
\begin{equation}
\|\nu(\cdot, r_n) - I_N\nu^n\|_1 + \eps^{2p}\|\partial_r \nu(\cdot, r_n) - I_N\mu^n\| \lesssim h^{m} + \kappa^{2}/\eps^{2p},\ 0 \leq n \leq \frac{T}{\kappa}.
\end{equation}
\end{theorem}
\begin{remark} Under the assumption of Theorem \ref{thm:HOE}, direct  error analysis for time-splitting schemes \cite{BaoFS,LU} would lead to the error estimates as $\|\nu(\cdot, r_n) - I_N\nu^n\|_1 + \eps^{2p}\|\partial_r \nu(\cdot, r_n) - I_N\mu^n\| \lesssim h^{m} +  \frac{\kappa^2}{\eps^{4p}}$. Our results are improved in the sense that the error bound $\frac{\kappa^2}{\eps^{4p}}$ is now $\frac{\kappa^2}{\eps^{2p}}$.
\end{remark}
\begin{remark}
The proof of the improved error bounds for the oscillatory complex NKGE in Theorem \ref{thm:HOE} is similar to the long-time problem and we omit the details for brevity. We will provide an example in section 5 to confirm the improved error bounds for the oscillatory complex NKGE and to demonstrate that they are sharp.
\end{remark}


\section{Numerical results}
In this section, we present some numerical examples in 1D and 2D to validate our improved uniform error bounds on the time-splitting methods for the long-time dynamics of the NKGE with weak nonlinearity and the improved error bounds for the oscillatory complex NKGE.

\subsection{The long-time dynamics in 1D}
First, we test the long-time errors of the TSFP \eqref{eq:gpsi}--\eqref{eq:gufull} for the NKGE \eqref{eq:GWNE} in 1D with $p=2$ and real-valued initial data as
\begin{equation}
u_0(x) = \frac{3}{2+\cos^2(x)},\quad u_1(x) = \frac{3}{4+\cos^2(x)}, \quad x \in\Omega= (0, 2\pi).	
\end{equation}
The numerical `exact' solution is computed by the TSFP \eqref{eq:gpsi}--\eqref{eq:gufull} with a very fine mesh size $h_e = \pi/60$ and  time step $\tau_e = 10^{-4}$. To quantify the error, we introduce the following error functions:
\begin{equation}
e_1(t_n) = \left\|u(x, t_n) - I_N u^n\right\|_1,	 \quad e_{1, \max}(t_n) = \max_{0 \leq q\leq n}e_1(t_q).
\end{equation}
In the rest of the paper, the spatial mesh size is always chosen sufficiently small such that the spatial errors can be neglected when considering the long-time temporal errors.

\begin{figure}[ht!]
\centerline{\includegraphics[width=14cm,height=5.5cm]{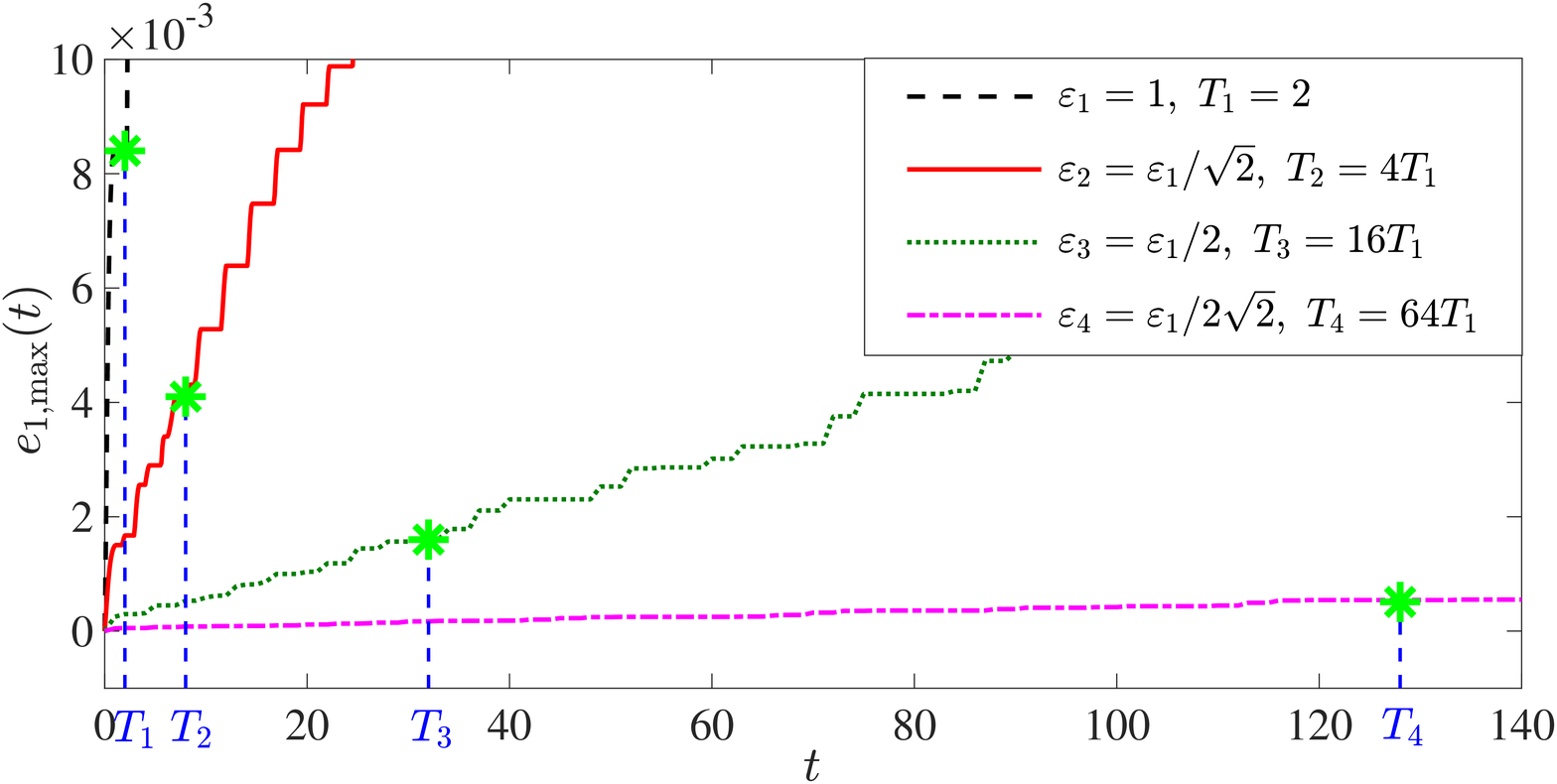}}
\caption{Long-time temporal errors of the TSFP \eqref{eq:gpsi}--\eqref{eq:gufull} for the NKGE \eqref{eq:GWNE} with $p=2$ and different $\varepsilon$ in 1D.}
\label{fig:1D_long}
\end{figure}

\begin{figure}[ht!]
\begin{minipage}{0.49\textwidth}
\centerline{\includegraphics[width=6.5cm,height=5.5cm]{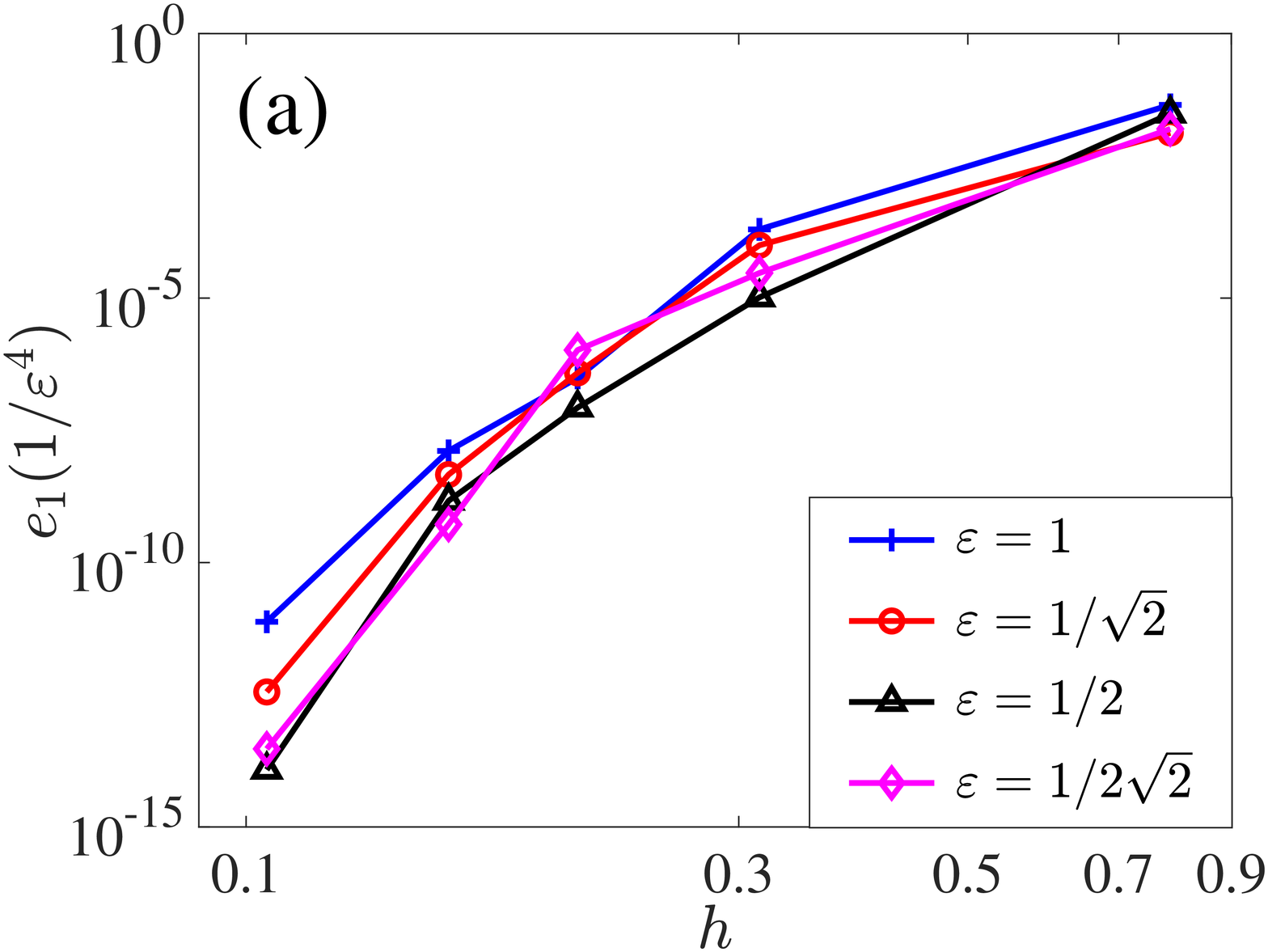}}
\end{minipage}
\begin{minipage}{0.49\textwidth}
\centerline{\includegraphics[width=6.5cm,height=5.5cm]{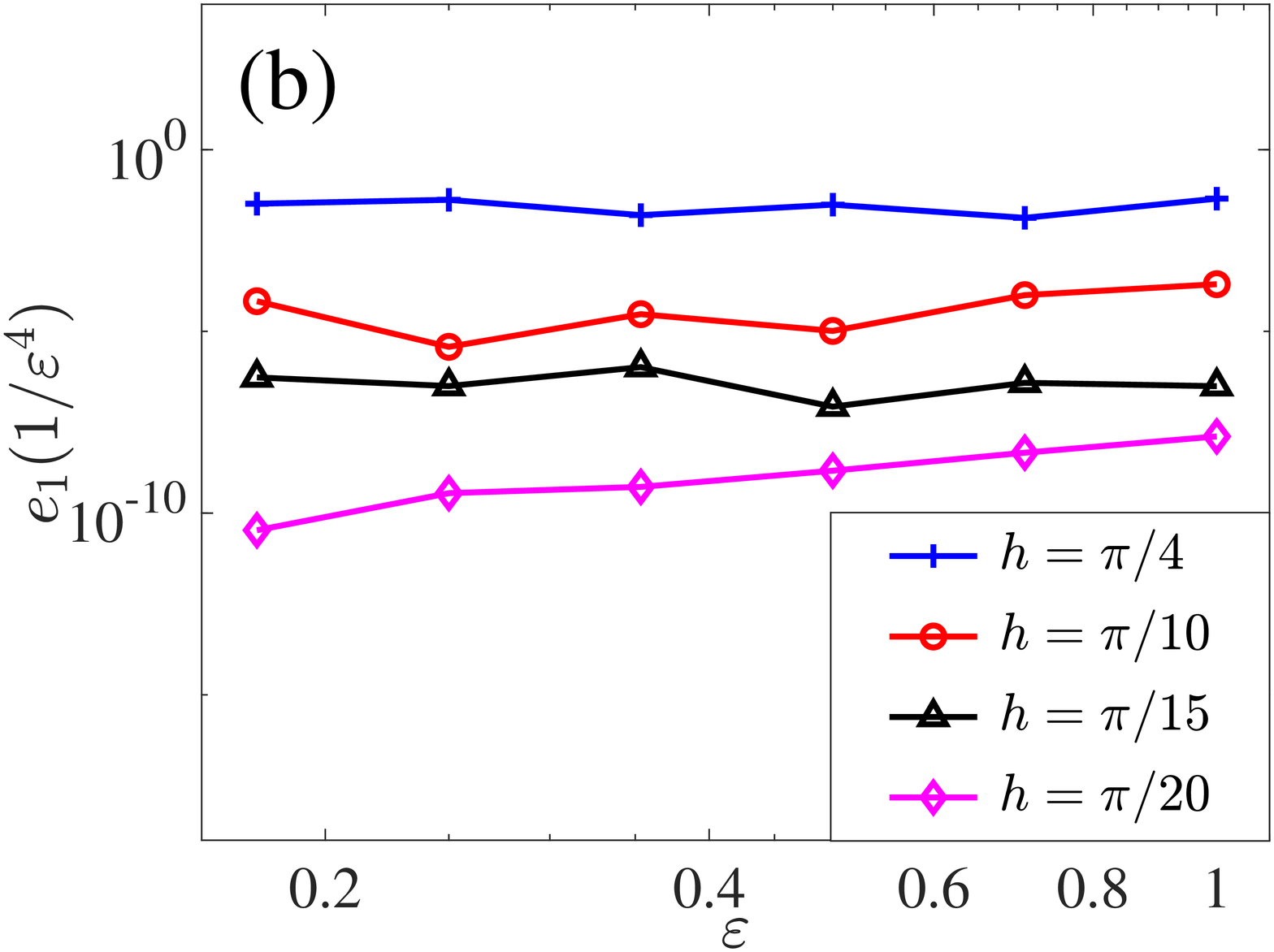}}
\end{minipage}
\caption{Long-time spatial errors of the TSFP \eqref{eq:gpsi}--\eqref{eq:gufull} for the NKGE \eqref{eq:GWNE} with $p = 2$ in 1D at $t = 1/\eps^4$.}
\label{fig:1D_spatial}
\end{figure}

\begin{figure}[ht!]
\begin{minipage}{0.49\textwidth}
\centerline{\includegraphics[width=6.5cm,height=5.5cm]{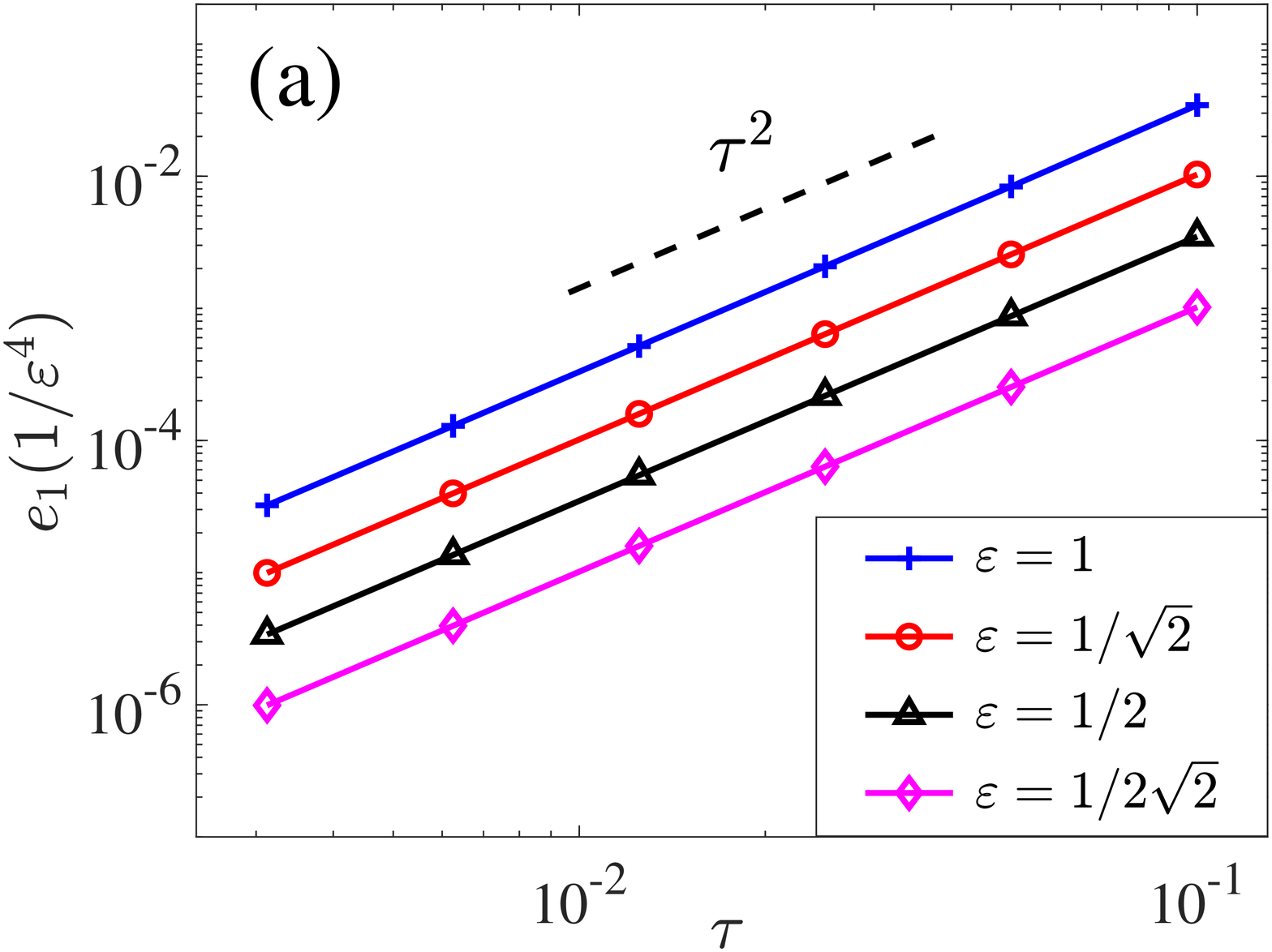}}
\end{minipage}
\begin{minipage}{0.49\textwidth}
\centerline{\includegraphics[width=6.5cm,height=5.5cm]{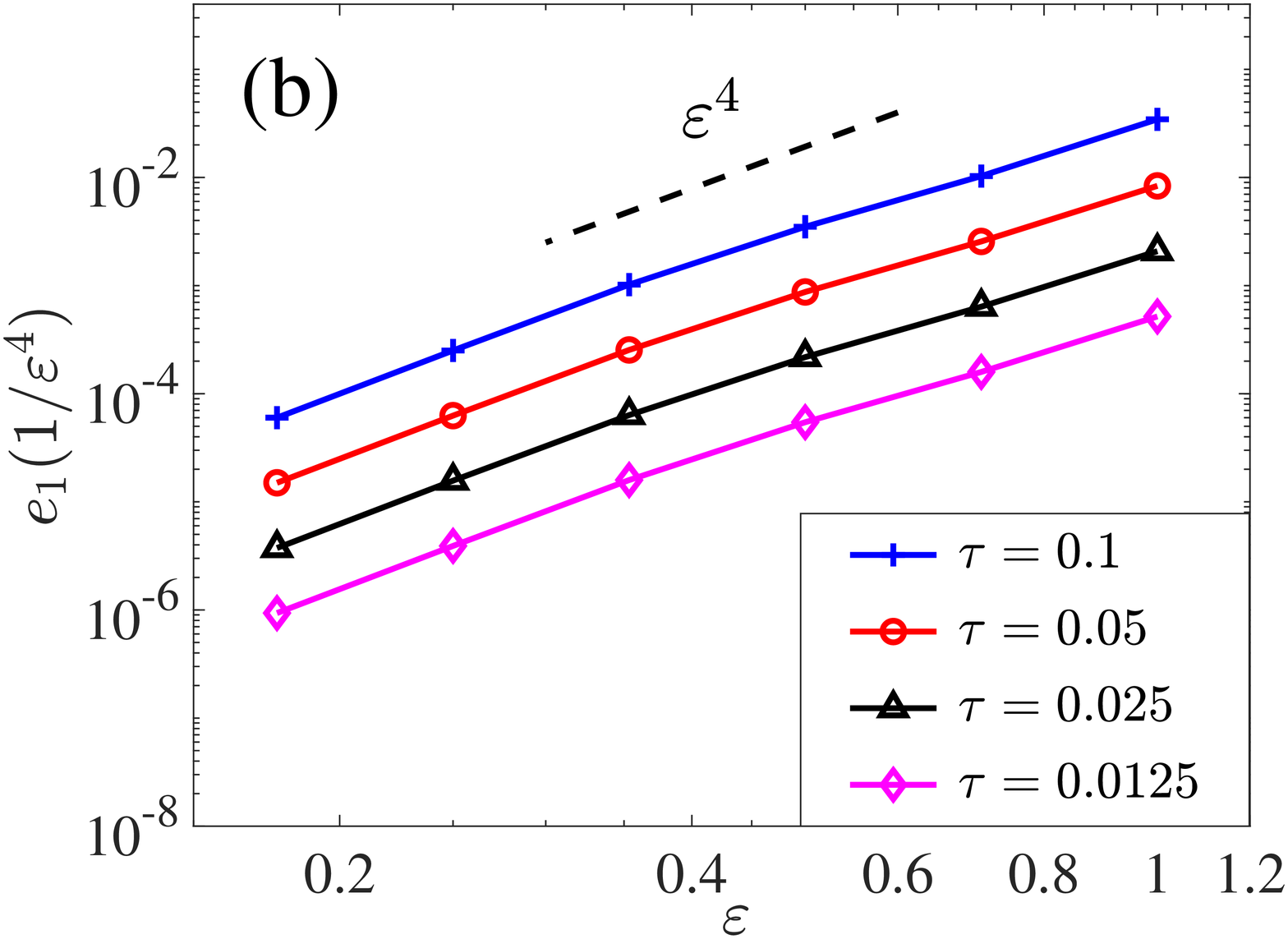}}
\end{minipage}
\caption{Long-time temporal errors of the TSFP \eqref{eq:gpsi}--\eqref{eq:gufull} for the NKGE \eqref{eq:GWNE} with $p = 2$ in 1D at $t = 1/\eps^4$.}
\label{fig:1D_temporal}
\end{figure}

Fig. \ref{fig:1D_long} displays the long-time errors of the TSFP \eqref{eq:gpsi}--\eqref{eq:gufull} for the NKGE \eqref{eq:GWNE} with $p =2$, the fixed time step $\tau$ and different $\varepsilon$, which confirms the improved uniform error bounds in $H^1$-norm at $O(\eps^4\tau^2)$ up to time at $O(1/\eps^4)$. Fig. \ref{fig:1D_spatial} and Fig. \ref{fig:1D_temporal} depict the spatial and temporal errors of the TSFP \eqref{eq:gpsi}--\eqref{eq:gufull} for the NKGE \eqref{eq:GWNE} with $p=2$ at $t = 1/\varepsilon^4$, respectively. Fig. \ref{fig:1D_spatial} indicates the spectral accuracy of the TSFP \eqref{eq:gpsi}--\eqref{eq:gufull} for the NKGE \eqref{eq:GWNE} in space and the spatial errors are independent of the small parameter $\varepsilon$. Each line in Fig. \ref{fig:1D_temporal} (a) corresponds to a fixed $\varepsilon$ and shows the global errors in $H^1$-norm versus the time step $\tau$, which confirms the second-order convergence of the TSFP \eqref{eq:gpsi}--\eqref{eq:gufull} for the NKGE \eqref{eq:GWNE} in time. Fig. \ref{fig:1D_temporal} (b) again validates that the global errors in $H^1$- norm behave like $O(\eps^4\tau^2)$ up to the time at $O(1/\eps^4)$.

For comparisons, we present the temporal errors of the first, second and fourth order splitting methods. In space, we use the Fourier pseudospectral method with a very fine mesh size such that the spatial errors are negligible.

\begin{figure}[ht!]
\begin{minipage}{0.49\textwidth}
\centerline{\includegraphics[width=6.5cm,height=5.5cm]{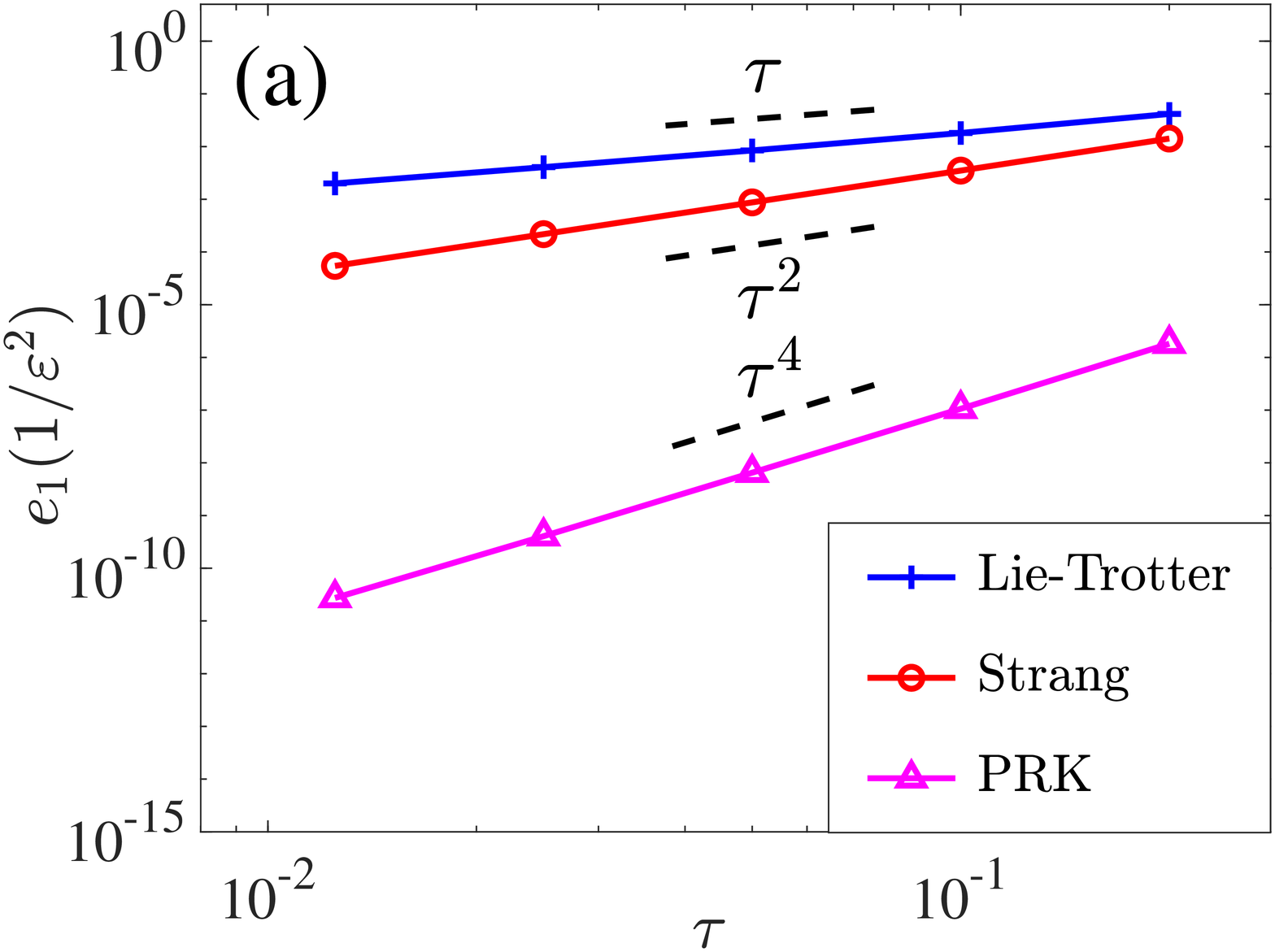}}
\end{minipage}
\begin{minipage}{0.49\textwidth}
\centerline{\includegraphics[width=6.5cm,height=5.5cm]{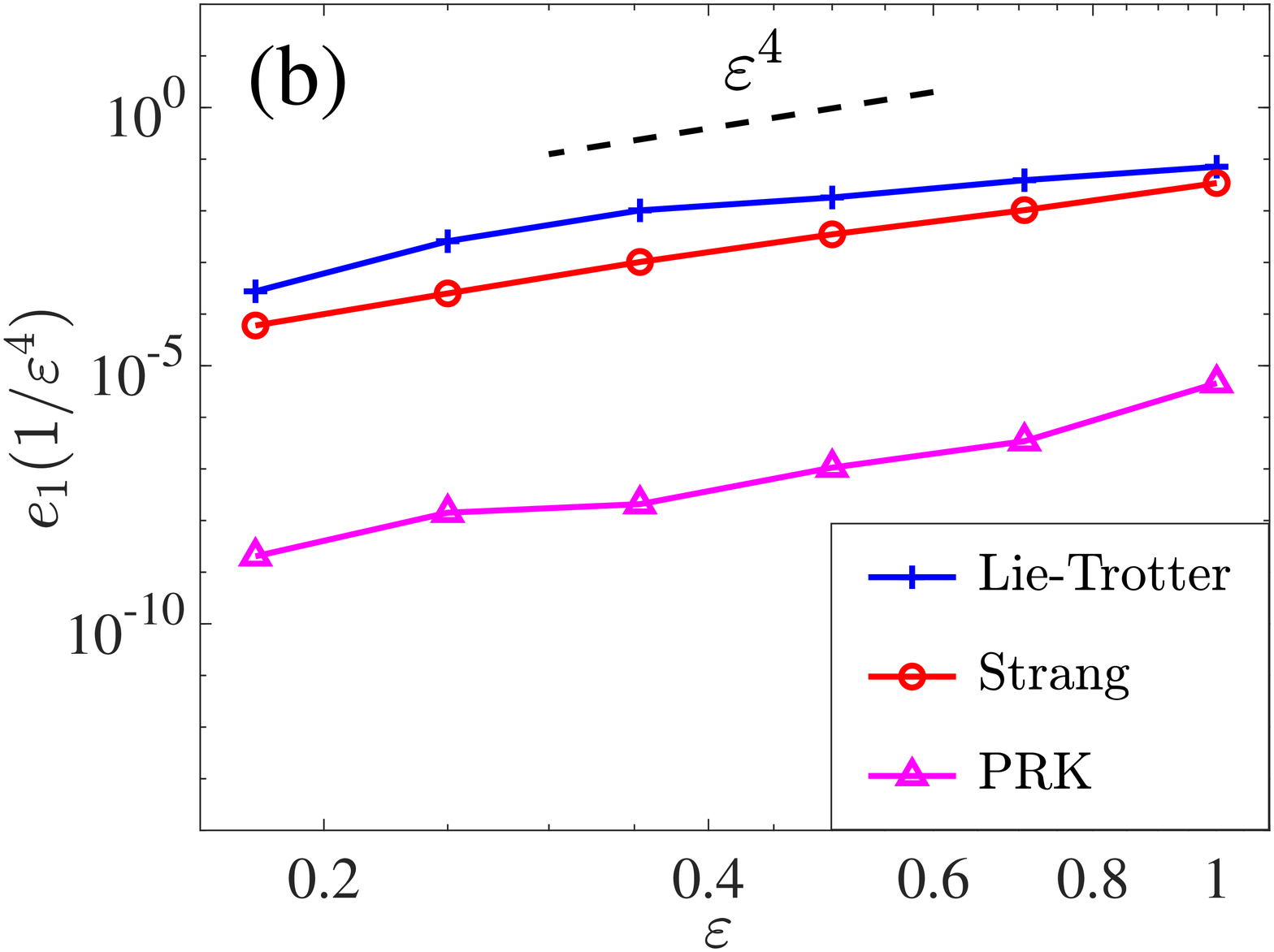}}
\end{minipage}
\caption{Comparisons of the first, second and fourth order splitting methods for the NKGE \eqref{eq:21}.}
\label{fig:Com}
\end{figure}

Fig. \ref{fig:Com} (a) depicts the temporal errors of three splitting methods with $\eps = 1/2$, which indicates that the higher order splitting method not only has higher order convergence rate but also achieves better accuracy under the same time step size. Fig. \ref{fig:Com} (b) shows the temporal errors of three splitting methods for the fixed time step and confirms the improved uniform error bounds for all the three splitting methods up to the time at $O(1/\eps^4)$.

\subsection{The long-time dynamics in 2D}
In this subsection, we show an example in 2D with the irrational aspect ratio of the domain $ (x, y) \in\Omega= (0, 1) \times (0, 2\pi)$. In the numerical experiment, we choose $p =1$ and the initial data as
\begin{equation*}
u_0(x, y) = \frac{2}{1 + \cos^2(2\pi x+y)},\quad u_1(x) = \frac{3}{2 + 2\cos^2(2\pi x+y)}.
\end{equation*}

\begin{figure}[ht!]
\centerline{\includegraphics[width=14cm,height=5.5cm]{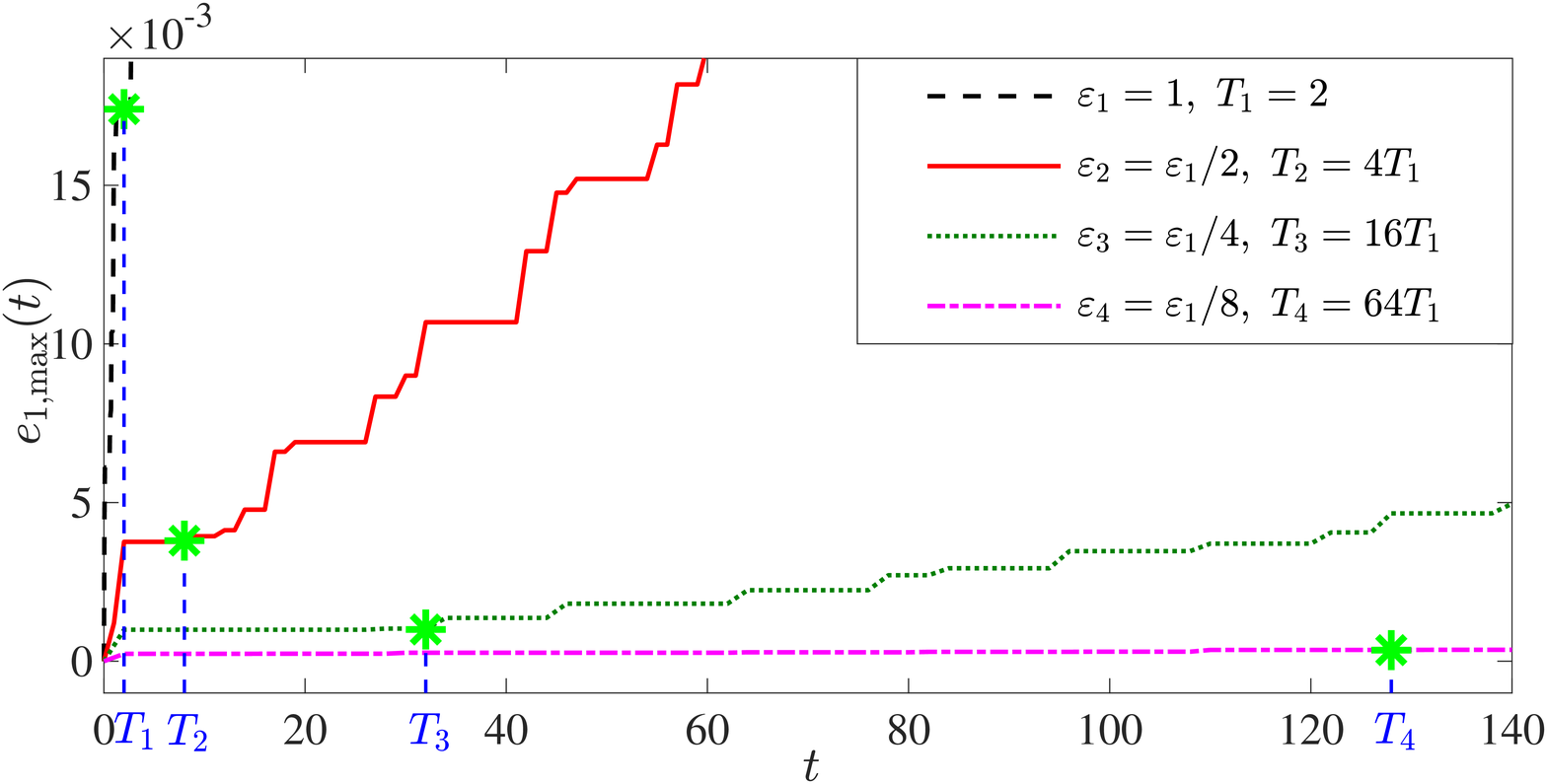}}
\caption{Long-time temporal errors of the TSFP method for the NKGE in 2D with different $\varepsilon$.}
\label{fig:2D_long}
\end{figure}

\begin{figure}[ht!]
\begin{minipage}{0.49\textwidth}
\centerline{\includegraphics[width=6.5cm,height=5.5cm]{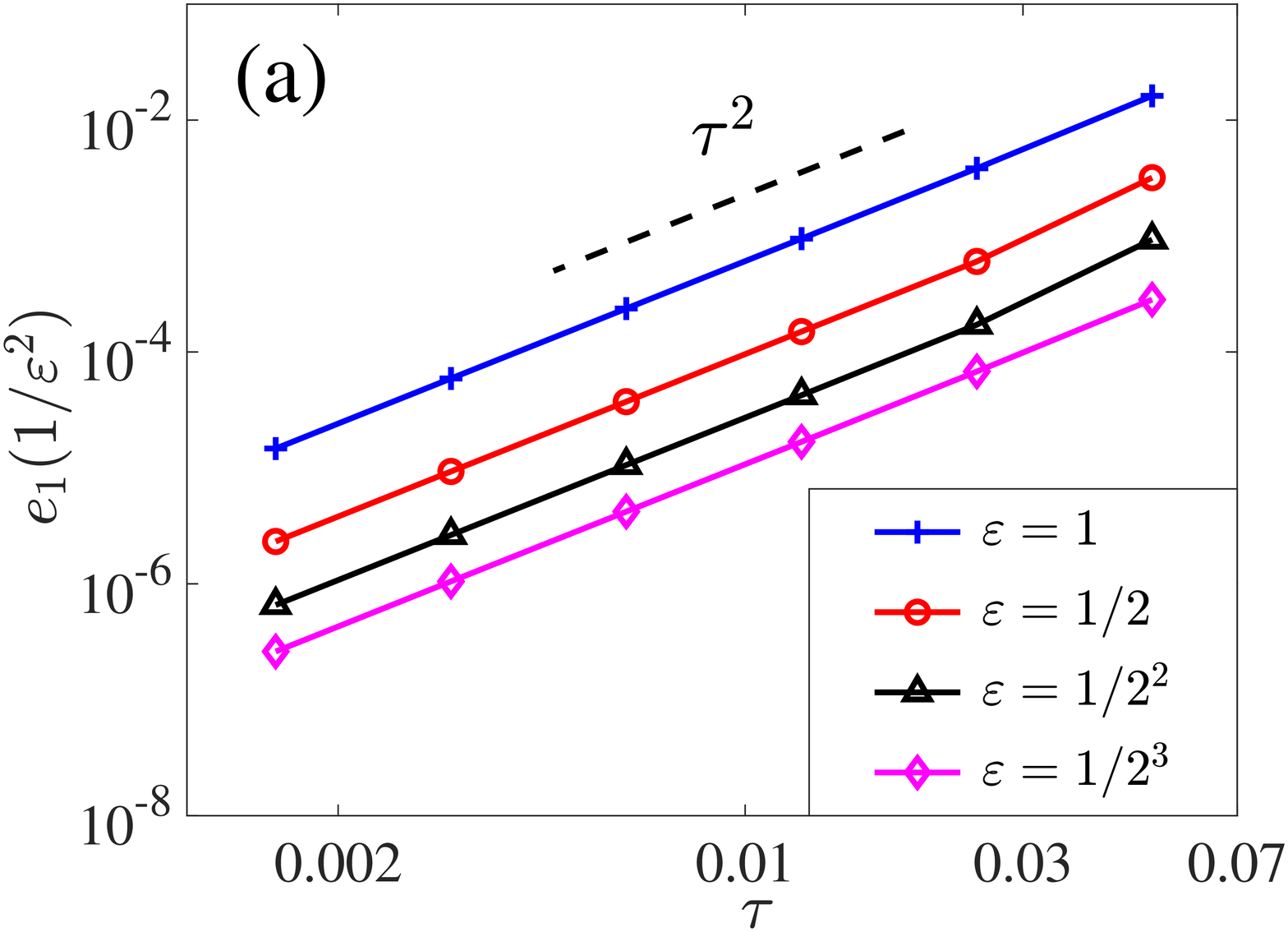}}
\end{minipage}
\begin{minipage}{0.49\textwidth}
\centerline{\includegraphics[width=6.5cm,height=5.5cm]{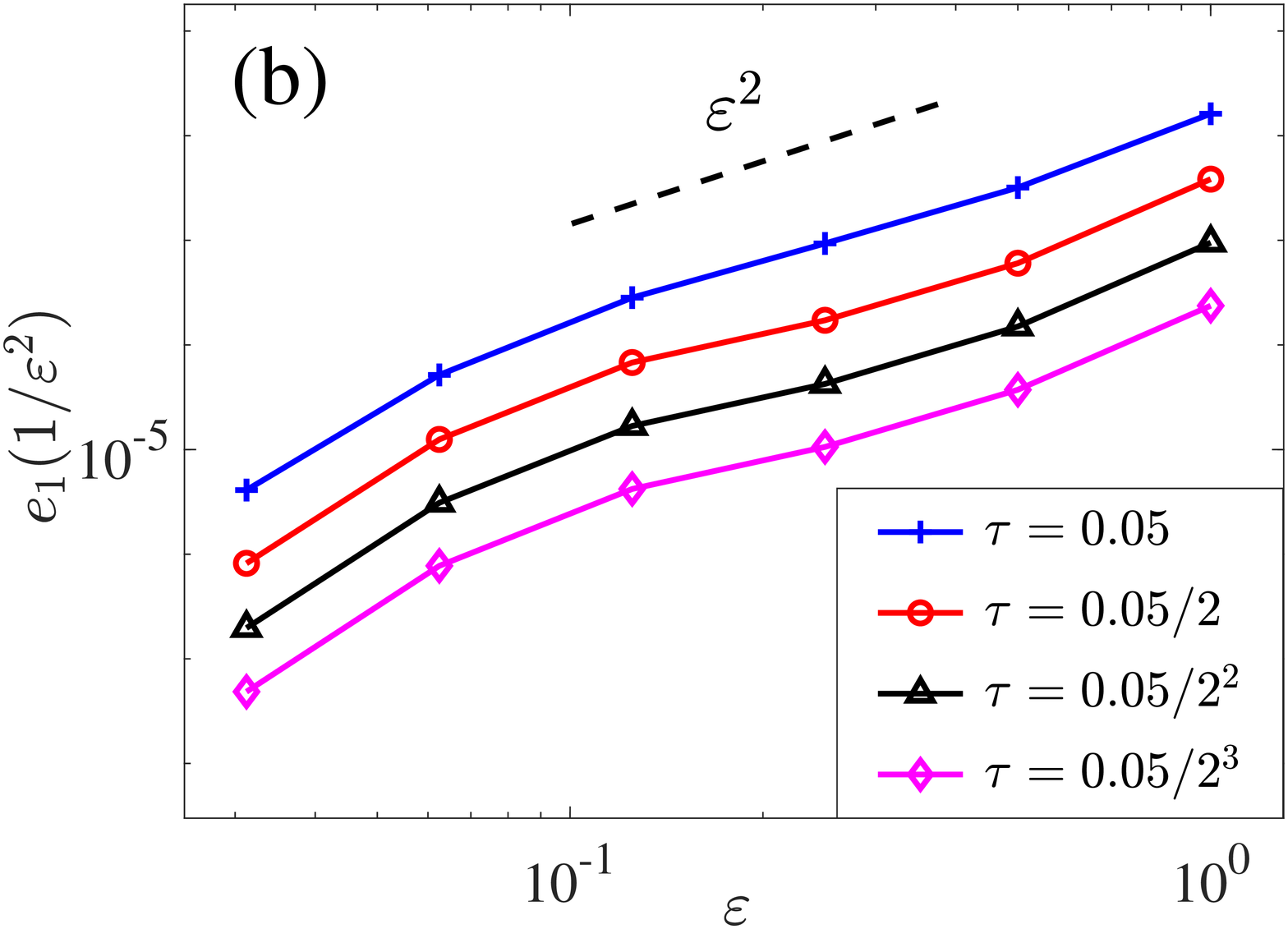}}
\end{minipage}
\caption{Long-time temporal errors of the TSFP method for the NKGE in 2D at $t = 1/\eps^2$.}
\label{fig:2D_temporal}
\end{figure}

Fig. \ref{fig:2D_long} presents the long-time errors of the TSFP method for the NKGE in 2D with a fixed time step $\tau$ and different $\varepsilon$, which confirms that the improved uniform error bounds at $O(\varepsilon^2\tau^2)$ up to the time at $O(1/\varepsilon^2)$ are also suitable for the domain with  irrational aspect ratio. Fig. \ref{fig:2D_temporal} depicts the temporal errors for the TSFP method for the NKGE in 2D at $t = 1/\varepsilon^2$, which again indicates that the TSFP method is second-order in time and validates the improved uniform error bounds up to the time at $O(1/\varepsilon^2)$.

\subsection{The oscillatory complex NKGE}
In this subsection, we present the numerical result for the oscillatory complex NKGE \eqref{eq:HOE_error} in 1D to confirm the improved error bound \eqref{eq:HOE}. We choose $p = 1$ and the complex-valued initial data as
\begin{equation*}
u_0(x) = x^2(x-1)^2+3,\quad u_1(x) = x(x-1)(2x-1) + 3i\cos(2\pi x),\quad x \in\Omega= (0, 1). 	
\end{equation*}
 The regularity is enough to ensure the improved error bound in $H^1$-norm.

\begin{table}
\caption{Temporal errors of the TSFP method for the oscillatory complex NKGE \eqref{eq:HOE} in 1D.}
\centering
\renewcommand\arraystretch{1.2}
\begin{tabular}{cccccc}
\hline
$e_1(r= 1)$ &$\kappa_0 = 0.05 $ & $\kappa_0/4 $ &$\kappa_0/4^2 $ & $\kappa_0/4^3 $ & $\kappa_0/4^4$   \\
\hline
$\varepsilon_0 = 1$ & \bf{1.11E-2} & 6.90E-4 & 4.31E-5 & 2.69E-6 & 1.68E-7  \\
order & \bf{-} & 2.00 & 2.00 & 2.00 & 2.00 \\
\hline
$\varepsilon_0 / 2 $  & 6.25E-2 & \bf{3.45E-3} & 2.14E-4 & 1.34E-5 & 8.35E-7  \\
order & -  & \bf{2.09} & 2.01 & 2.00 & 2.00 \\
\hline
$\varepsilon_0 / 2^2 $ & 8.26E-1 & 1.89E-2 & \bf{1.11E-3} & 6.93E-5 & 4.33E-6 \\
order & -  & 2.72 & \bf{2.04} & 2.00 & 2.00 \\
\hline
$\varepsilon_0 / 2^3 $ & 1.54 & 3.19E-1 & 1.62E-2 & \bf{1.01E-3} & 6.29E-5 \\
order & -  & 1.14 & 2.15  & \bf{2.00} & 2.00 \\
\hline
$\varepsilon_0 / 2^4$ & 2.09 & 3.82 & 7.90E-2 & 4.26E-3 & \bf{2.64E-4}  \\
order & - & -0.44 & 2.80 & 2.11 & \bf{2.01} \\
\hline
\end{tabular}
\label{tab:HOE}
\end{table}

Table \ref{tab:HOE} lists the temporal errors of the TSFP method for the oscillatory NKGE \eqref{eq:HOE} in 1D, which indicates that the second-order convergence can only be observed when $\kappa\lesssim \eps^2$ (cf. the upper triangle above the diagonal with bold letters) and the temporal errors in $H^1$-norm behave like $O(\kappa^2/\eps^2)$ to confirm the improved error bound \eqref{eq:HOE_error} and to demonstrate that they are sharp.

\section{Conclusions}
Improved uniform error bounds on the time-splitting methods for the long-time dynamics of the nonlinear Klein--Gordon equation (NKGE) with weak cubic nonlinearity were rigorously established. By employing the technique of regularity compensation oscillation (RCO), the improved uniform error bounds for the second-order semi-discretization and full-discretization up to the time at $O(1/\varepsilon^2)$ were carried out at $O(\varepsilon^2\tau^2)$ and  $O(h^m + \eps^2\tau^2)$, respectively. The improved error bounds are extended to the complex NKGE with a general power nonlinearity in the long-time regime and the oscillatory complex NKGE up to the fixed time $T$. Numerical results in 1D and 2D were presented to confirm the improved error bounds and to demonstrate that they are sharp.


\bibliographystyle{siamplain}

\end{document}